\documentclass[12pt]{article}
\title{\Large Structural Properties of ${\cal R}_2$ Part II}
\author{\normalsize Timothy J. Carlson \\
  \normalsize The Ohio State University, Columbus, OH 43210 USA \\
 \normalsize email: carlson@math.ohio-state.edu 
}

\date{\vspace{-0.3in}}

\vspace{-.2in}
\usepackage{amsmath}    
\usepackage{amssymb}
\usepackage{graphicx}   
\usepackage{verbatim}   
\usepackage{color}      
\usepackage{subfigure}  
\usepackage{hyperref}   

\newcommand{\qed}{\hfill {$\square$}}
\newcommand{\blankline}{\vspace{4 mm}}
\newcommand{\halfblankline}{\vspace{2 mm}}

\newcommand{\calRrho}{\mbox{${\cal R}^\rho$}}
\newcommand{\calRtwo}{\mbox{${\cal R}_2$}}
\newcommand{\calRtwoalpha}{\mbox{${\cal R}_2^\alpha$}}

\newcommand{\leqonerho}{\mbox{$\leq_1^\rho$}}

\newcommand{\leqtwo}{\mbox{$\leq_2$}}

\newtheorem{prop}{Proposition}[section]

\newtheorem{thm}[prop]{Theorem}
\newtheorem{lem}[prop]{Lemma}

\newtheorem{dfn}[prop]{Definition}

\begin{document}

\maketitle

\blankline

\noindent
{\small 
{\bf Abstract.} 
This is the second of two papers establishing structural properties of \calRtwo, the structure giving rise to pure patterns of resemblance of order two, which partially underly the results in [\ref{Ca09}] and [\ref{Ca16}] as well as other work in the area. 
}

\blankline

For the entire paper, we will assume $\rho$ is an arbitrary additively indecomposable ordinal and $\alpha$ is an epsilon number greater than $\rho$.
By Lemma 5.8 of [\ref{Ca??}], $\kappa^\rho_\alpha$ is additively indecomposable.

The {\it Order Reduction Theorem} says that for $\eta<\theta^\rho_2(\alpha)$, $J^\rho_{\alpha,\eta}$ is essentially isomorphic to an initial segment of ${\cal R}^\alpha_2$.
The important point is that this reduction reduces the lengths of chains in \leqtwo: if the longest chain of elements of $J^\rho_{\alpha,\eta}$ with respect to \leqtwo\ has length $n+1$ where $n\in \omega$ then the longest chain among the elements of the corresponding initial segment of ${\cal R}^\alpha_2$ with respect to the analogue of \leqtwo\ in ${\cal R}^\alpha_2$ has length at most $n$.

The {\it Second Recurrence Theorem for \leqtwo} says roughly that the initial segment of ${\cal R}^\alpha_2$ corresponding to $I_{\omega^\eta}$ grows by adding $\rho$ new intervals of the form  $I^\alpha_\beta$ each time $\eta$ increases.
In the case $\rho=1$, which is essentially \calRtwo, this says one new interval is added each time $\eta$ increases.

Section \ref{sp2} contains background.

Section \ref{sli} introduces the notion of {\it local incompressibility}.

Section \ref{sor} establishes the Order Reduction Theorem.

Section \ref{ssecondrtwo} establishes the Second Recurrence Theorem for \leqtwo.

\section{Preliminaries}
\label{sp2}

We will use the notation   from [{\ref{Ca??}] except we will need to be more careful about parameters in definitions as we now discuss.
As in [\ref{Ca??}], our metatheory is ${\sf KP}\omega$, Kripke-Platek set theory
with the axiom of infinity (see [\ref{Ba75}]).

Assume $X$ is a set of ordinals and $\sigma$ is an additively indecomposable ordinal. 
We will say that $X$ is  {\it $\sigma$-closed} if it is closed in ${\cal R}^\sigma$ (recall the definition of {\it closed} from [\ref{Ca??}] which is not the same as being topologically closed).
We will say that a map $h:X\rightarrow ORD$ is a {\it $\sigma$-covering} of $X$ if $h$ is a covering of $X$, as a substructure of ${\cal R}_2^\sigma$, in ${\cal R}^\sigma_2$. 
Similarly, $\tilde{X}$ is a {\it $\sigma$-covering} of $X$ if $\tilde{X}$ is a covering of $X$ when $X$ and $\tilde{X}$ are viewed as a substructures of ${\cal R}_2^\sigma$.
Similary, we will say that $h$ is a {\it $\sigma$-embedding} of $X$ if $h$ is an embedding of $X$, as a substructure of ${\cal R}^\sigma_2$, into ${\cal R}^\sigma_2$.
We will say that $h$ is a {\it $\sigma$-isomorphism} of $X$ and $Y$ if $h$ is an isomorphism of $X$ and $Y$ as substructures of ${\cal R}^\sigma_2$.
We will say $X$ and $Y$ are {\it $\sigma$-isomorphic} if there is a $\sigma$-isomorphism of $X$ and $Y$.
We will also write $X\cong_\sigma Y$ when $X$ and $Y$ are $\sigma$-isomorphic.
On the other hand, we will almost always use the notation $\kappa_\xi^\sigma$, $I_\xi^\sigma$, $\theta_1^\sigma$, $index^\sigma$, $\theta_2^\sigma(\beta)$, $\nu^\sigma_{\beta,\xi}$, $J^\sigma_{\beta,\xi}$ and $\overline{J}^\sigma_{\beta,\xi}$ in the case $\sigma=\rho$ and $\beta=\alpha$. 
So, we continue to suppress parameters in this case and write $\kappa_\xi$, $I_\xi$, $\theta_1$, $index$, $\theta_2$, $\nu_\xi$, $J_\xi$ and $\overline{J}_\xi$ respectively.

\section{Local Incompressibility}
\label{sli}


Recall that the Recurrence Theorem for Small Intervals (Theorem 8.1 of [\ref{Ca??}]), which will be abbreviated as RTSI, implies the following.
\begin{itemize}
\item
For all $\delta$ and $0<\chi_1\leq\chi_2<\kappa_\alpha$ 
\begin{center}
$\kappa_\alpha\cdot\delta+\chi_1\leq^\rho_k \kappa_\alpha\cdot\delta+\chi_2$ \ \ \ iff \ \ \ $\chi_1\leq^\rho_k \chi_2$
\end{center}
\item 
For all $\delta$ and all $0<\chi<\kappa_\alpha$
 \begin{itemize}
 \item
 $\kappa\not\leq_2^\rho \kappa_\alpha\cdot\delta +\chi$ whenever $\kappa\leq\kappa_\alpha\cdot\delta$.
 \item
 $\kappa_\alpha\cdot\delta + \chi \not\leq^\rho_1 \kappa_\alpha\cdot\delta +\kappa_\alpha$.
 \end{itemize}
 \item
 For all ordinals $\delta_1$,  $\delta_2$ and $\chi$ with $\delta_1\leq\delta_2$ and $0<\chi<\kappa_\alpha$,
 \begin{center}
 $\kappa_\alpha\cdot\delta_1\leq^\rho_1\kappa_\alpha\cdot\delta_2+\chi$ \ \ iff \ \ $\kappa_\alpha\cdot\delta_1\leq^\rho_1 \kappa_\alpha\cdot\delta_2 +\kappa_\gamma$
 \end{center}
 where $\gamma=index(\chi)$.
\end{itemize}
For the third part, use  the fact that $\kappa_\alpha\cdot\delta_2+\kappa_\gamma\leq^\rho_1 \kappa_\alpha\cdot\delta_2+\chi$ (by the first part).


\begin{dfn}
Assume $X$ is a set of ordinals. Let $M$ be the collection of ordinals $\mu$ such that $\mu$ is divisible by $\kappa_\alpha$ and $X\cap [\mu, \mu + \kappa_\alpha)\not=\emptyset$. 
For $\mu \in M$, there is a unique nonempty subset $R_\mu$ of $\kappa_\alpha$ such that $X\cap [\mu, \mu+ \kappa_\alpha)=\mu +R_\mu$. 
The indexed family $R_\mu$ $(\mu\in M)$ will be called the {\bf interval decomposition} of $X$.
If $0\in R_\mu$ and $R_\mu$ is incompressible for each $\mu\in M$, we will say that $X$ is {\bf locally incompressible}.
\end{dfn}

Notice that $X$ is $\rho$-closed iff $R_\mu$ is $\rho$-closed for all $\mu\in M$.
Moreover, $X$ is $\kappa_\alpha$-closed  iff $0\in R_\mu$ for all $\mu\in M$. 
Hence, if $X$ is locally incompressible then $X$ is both $\rho$-closed and $\kappa_\alpha$-closed.


\begin{lem} 
\label{lsrtt}
Assume $X$ is a finite $\rho$-closed set of ordinals with interval decomposition $R_\mu$ $(\mu\in M)$.
\begin{enumerate}
\item
Assume $X'$ has interval decomposition $R'_\mu$ $(\mu\in M)$.
If $R'_\mu$ is either a $\rho$-incompressible $\rho$-covering of $R_\mu$ or $R'_\mu=R_\mu$ for $\mu\in M$ then $X'$ is a $\rho$-covering of $X$.
\item
Assume $\mu_1$ and $\mu_2$ are divisible by $\kappa_\alpha$ and $\mu_1\leq\mu_2$.
If $\mu_1\leq X < \mu_2+\kappa_\gamma$ where $\gamma<\alpha$ then there is a locally incompressible set $X'$  such that $X\subseteq X'$ and $\mu_1\leq X'<\mu_2+\kappa_\gamma$.
\item
Assume $X$ is locally incompressible. 
If $h$ is a $\rho$-covering of $X$ such that $h(\mu)$ is divisible by $\kappa_\alpha$ for all $\mu\in M$ then 
$$h(\mu +\chi) \geq h(\mu)+\kappa_\gamma$$
whenever $\chi\in R_\mu$ and $\chi\in I_\gamma$.
\end{enumerate}
\end{lem}
{\bf Proof.}
For part 1, let $h_\mu$ be the $\rho$-covering of $R_\mu$ onto $R'_\mu$.
Define $h$ on $X$ by $h(\mu +\chi)=\mu +h_\mu(\chi)$ whenever $\chi\in R_\mu$.
Since $\rho$ divides $\kappa_\alpha$, Lemma 2.4 of [\ref{Ca??}] implies $h$ is an embedding of $X$ into \calRrho. 

Notice that if $\mu\in M\cap X$ then $h(\mu)=\mu+h_\mu(0)=\mu+0=\mu$ (since $index(h_\mu(0))\leq index(0)=0$).

To show $h$ is a $\rho$-covering, assume $\chi_i\in R_{\mu_i}$ for $i=1,2$ and $\mu_1+\chi_1<^\rho_k \mu_2+\chi_2$.
We will show $h(\mu_1+\chi_1)\leq^\rho_k h(\mu_2+\chi_2)$.

First consider the case when $\chi_1\not=0$.
By RTSI, $\mu_1=\mu_2$.
Let $\mu$ be the common value.
By RTSI again, $\chi_1\leq^\rho_k \chi_2$.
Since $h_\mu$ is a $\rho$-covering, $h_\mu(\chi_1)\leq^\rho_k h_\mu(\chi_2)$.
By RTSI once again, $\mu+h_\mu(\chi_1)\leq^\rho_k \mu+h_\mu(\chi_2)$ i.e. $h(\mu_1+\chi_1)\leq^\rho_k h(\mu_2+\chi_2)$.

Now consider the case when $\chi_1=0$. If $\chi_2=0$ then $h(\mu_i)=\mu_i$ for $i=1,2$ from which the desired conclusion is trivial.
Assume $\chi_2\not=0$.
By RTSI, $k=1$.
Let $\gamma=index(\chi_2)$.
So, $\kappa_\gamma\leq \chi_2 \leq max_1(\kappa_\gamma)$.
By RTSI, $\mu_2+\kappa_\gamma \leqonerho \mu_2+max_1(\kappa_\gamma)$.
Since $\mu_1\leqonerho \mu_2+\chi_2$, $\mu_1\leqonerho \mu_2 + max_1(\kappa_\gamma)$.
Since $index(h_{\mu_2}(\chi_2))\leq \gamma$, $h_{\mu_2}(\chi_2)\leq max_1(\kappa_\gamma)$.
Therefore, $h(\mu_1)=\mu_1 \leqonerho \mu_2 + h_{\mu_2}(\chi_2) =h(\mu_2+\chi_2)$.

For part 2,  let $R_\mu$ $(\mu\in M)$ be the interval decomposition of $X$.
For each $\mu\in M$, let $R'_\mu$ be a $\rho$-incompressible set such that $0\in R'_\mu$, $R_\mu\subseteq R'_\mu$ and $max(index[R_\mu])=max(index[R'_\mu])$.
Let $X'=\bigcup_{\mu\in M}R'_\mu$.
One easily checks $X'$ has the desired properties.

For part 3, assume $\mu\in M$.
Let $Q_1=\{\chi\in R_\mu \, | \, h(\mu+\chi)<h(\mu)+\kappa_\alpha\}$ and $Q_2= R_\mu - Q_1$.
Clearly, $Q_1<Q_2$. 

We claim $Q_1\not\leq^\rho_1 Q_2$.
Assume $\chi_i \in Q_i$ for $i=1,2$.
We will show $\chi_1\not\leq^\rho \chi_2$.
Since this is clear if $\chi_1=0$, we may assume $0<\chi_1$.
Therefore, 
$$h(\mu)<h(\mu+\chi_1)<h(\mu)+\kappa_\alpha \leq h(\mu+\chi_2)$$
By RTSI, $h(\mu +\chi_1)\not\leq^\rho_1 h(\mu + \chi_2)$.
Since $h$ is a $\rho$-covering, $\mu+\chi_1\not\leq^\rho_1 \mu +\chi_2$.
By RTSI again, $\chi_1\not\leq^\rho_1 \chi_2$.

By part 8 of Lemma 6.4 of [\ref{Ca??}], $Q_1$ is $\rho$-incompressible.
Define $f:Q_1\rightarrow \kappa_\alpha$ by 
$$h(\mu+\chi)=h(\mu)+f(\chi)$$ 
We claim $f$ is a $\rho$-covering of $Q_1$.
Assume $\chi_1,\chi_2\in Q_1$ and $\chi_1<^\rho_k \chi_2$.
The assumption implies $0<\chi_1$.
By RTSI, $\mu +\chi_1 \leq^\rho_k \mu+\chi_2$. 
Since $h$ is a $\rho$-covering, RTSI implies that $f(\chi_1)\leq^\rho_k f(\chi_2)$.

To prove the conclusion of part 2, assume $\chi\in R_\mu$ and $\chi\in I_\gamma$. 
If $\chi\in Q_2$ the conclusion of part 2 is clear.
So, we may assume $\chi\in Q_1$.
Since $f$ is a $\rho$-covering of $Q_1$ and $Q_1$ is $\rho$-incompressible, $index(f(\chi))\geq index(\chi)=\gamma$.
Therefore, $f(\chi)\geq \kappa_\gamma$.
By definition of $f$, this implies the desired conclusion.
\qed


\begin{lem} 
\label{lorA0}
Assume $\lambda$ is a limit ordinal, $\nu<\lambda$ is divisible by $\kappa_\alpha$, $X\subseteq \nu$ is a finite $\rho$-closed set and $\cal Y$ is a family of nonempty finite $\rho$-closed subsets of $\lambda$  such that 
\begin{itemize}
\item
$\cal Y$ is cofinal in $\lambda$.
\item
The elements of $\cal Y$ have the same cardinality.
\item 
$rem^\rho[Y_1]=rem^\rho[Y_2]$ for all $Y_1,Y_2\in \cal Y$.
\end{itemize}
For $Y\in \cal Y$, let $M_Y$ be the set of $\mu$ which are divisible by $\kappa_\alpha$ such that $Y$ intersects $[\mu,\mu+\kappa_\alpha)$.
There exists finite $R\subseteq \kappa_\alpha$ and ${\cal Y}'\subseteq \cal Y$ such that
\begin{enumerate}
\item
$R$ is $\rho$-incompressible and $0\in R$.
\item
${\cal Y}'$ is cofinal in $\lambda$.
\item
There is a $\rho$-covering of $X\cup Y$ into $X\cup (M_Y+R)$ for all $Y\in {\cal Y}'$.
\end{enumerate}
\end{lem}
{\bf Proof.}
Assume $Y\in \cal Y$. 
Let $R_{Y,\mu}$ $(\mu\in M_Y)$ be the interval decomposition of $Y$.
Let $R_Y=\bigcup_{\mu\in M_Y}R_{Y,\mu}\cup \{0\}$.
Clearly, $Y\subseteq M_Y+R_Y$. 
Since the elements of $M_Y$ are divisible by $\kappa_\alpha$ which in turn is divisible by $\rho$,
$rem^\rho[Y]=\bigcup_{\mu\in M_Y}rem^\rho[M_{Y,\mu}]$.
Therefore, $rem^\rho[R_Y]=\bigcup_{\mu\in M_Y} rem^\rho[M_{Y ,\mu}] \cup \{0\} = rem^\rho[Y]\cup \{0\}$.
Since $rem^\rho[Y]$ is fixed as $Y$ varies over $\cal Y$, so is $rem^\rho[R_Y]$. 

Lemma 3.3 of [\ref{Ca??}] implies there exists ${\cal Y}'\subseteq \cal Y$ which is cofinal in $\lambda$ such that $R_{Y_1}\cong_\rho R_{Y_2}$ for all $Y_1,Y_2\in {\cal Y}'$.
Moreover, we may assume $\nu\leq Y$ for $Y\in {\cal Y}'$.
This implies that $X<M_Y$ for $Y\in {\cal Y}'$.
Let $R$ be a $\rho$-incompressible covering of $R_Y$ for $Y\in {\cal Y}'$.
Since $R_Y$ is a subset of $\kappa_\alpha$, so is $R$. 
Since $0\in R_Y$, $0\in R$.

Assume $Y\in {\cal Y}'$.
Part 1 of Lemma \ref{lsrtt} implies $X\cup (M_Y+R)$ is a covering of $X\cup (M_Y+R_Y)$.
This implies there is a covering of $X\cup Y$ into $X\cup (M_Y+R)$.
\qed

\section{The Order Reduction Theorem}
\label{sor}

Recall that we are assuming $\rho$ is an arbitrary additively indecomposable ordinal and $\alpha$ is an epsilon number greater than $\rho$.

We will see that if $\xi_1<\xi_2<\theta_2$ are additively indecomposable then $J_{\xi_1}$ is isomorphic to a proper initial segment of $J_{\xi_2}$.
For the proof, we collapse the \leqonerho\ connectivity components of each small interval to single points.
More precisely, in Definition \ref{drec2D} below the interval $\nu_\xi +\kappa_\alpha\cdot\delta+I_\gamma$ is replaced by the single point $\nu_\xi +\kappa_\alpha\cdot\delta+\kappa_\gamma$ when $\gamma<\alpha$ and $\nu_\xi +\kappa_\alpha\cdot\delta+\kappa_\gamma\in J_\xi$.
We will see that the resulting structure is essentially isomorphic to an initial segment of \calRtwoalpha\ and use that result to establish a second recurrence theorem for the length of $J_\xi$ in the following section.

\begin{dfn} 
\label{drec2D}
Define $\overline{J}'_\xi$ to be the collection of elements of $\overline{J}_\xi$ which are of the form $\nu_\xi+\kappa_\alpha\cdot\delta+\kappa_\gamma$ where $\gamma<\alpha$.
The operation $\iota_{\alpha,\xi}$  is the unique order preserving map of an initial segment of $ORD$ onto $\overline{J}'_\xi$. 
Define the function $\varphi_\xi:\overline{J}_\xi\rightarrow dom(\iota_\xi)$ by
$$\varphi_\xi(\nu_\xi+\kappa_\alpha\cdot\delta +\chi)\ = \ \alpha\cdot\delta + index(\chi)$$
when $\chi<\kappa_\alpha$.
\end{dfn}

As usual, we should include the parameter $\rho$ in the notation for $\iota_{\alpha,\xi}$ and write $\iota^\rho_{\alpha,\xi}$.
Similarly for $\varphi_\xi$.
Since both $\rho$ and $\alpha$ are clear from the context, we will simply write $\iota_\xi$ and $\varphi_\xi$.

The following observations will be used often.
\begin{itemize}
\item For all ordinals $\alpha\cdot\delta+\gamma$ in the domain of $\iota_\xi$ with $\gamma<\alpha$,
$$\iota_\xi(\alpha\cdot\delta+\gamma)=\nu_\xi+\kappa_\alpha\cdot\delta+\kappa_\gamma$$
implying
$$\iota_\xi(\alpha\cdot\delta+\gamma) = \iota_\xi(\alpha\cdot\delta)+ \kappa_\gamma$$
\item
If $\sigma\in dom(\iota_\xi)$ is divisible by $\alpha$ then $\iota_\xi(\sigma)$ is divisible by $\kappa_\alpha$.
\item
If $\mu\in \overline{J}_\xi$ is divisible by $\kappa_\alpha$ then $\mu\in \overline{J}_\xi'$ and $\iota_\xi^{-1}(\mu)$ is divisible by $\alpha$.
\item
$\iota_\xi$ is continuous and the domain of $\iota_\xi$ has a maximal element if it is bounded.
\item
$\varphi_\xi$ is weakly order preserving i.e. if $\sigma,\tau\in \overline{J}_\xi$ and $\sigma\leq \tau$ then $\varphi_\xi(\sigma)\leq \varphi_\xi(\tau)$.
\item
$\varphi_\xi$ extends $\iota_\xi^{-1}$ and
$$\varphi_\xi(\mu+\chi)=\iota_\xi^{-1}(\mu)+index(\chi)$$
whenever $\mu+\chi\in\overline{J}_\xi$, $\mu$ is divisible by $\kappa_\alpha$ and $\chi<\kappa_\alpha$.
\item
If $\lambda$ is in the domain of $\iota_\xi$ then $\iota_\xi(\lambda)$ is the least $\sigma$ such that $\varphi_\xi(\sigma)=\lambda$.
\item
If $Z$ is an $\alpha$-closed subset of the domain of $\iota_\xi$ then $\iota_\xi[Z]$ is $\kappa_\alpha$-closed.
\item
  If $Z$ is a $\kappa_\alpha$-closed subset of $ \overline{J}_\xi$ then
 $\varphi_\xi[Z]$ is $\alpha$-closed.
  \end{itemize}

The second observation follows from the continuity of $\gamma\mapsto \kappa_\gamma$.

The following lemma will be the key to the inductive proof of the Order Reduction Theorem.


\begin{lem} 
\label{lorA}
Assume $\xi<\theta_2$, $\lambda\in dom(\iota_\xi)$ and 
\begin{itemize}
\item[$(\star)$]
\hspace{1.1in} $\sigma\leq^\alpha_k \tau \ \ \ \ \Longleftrightarrow \ \ \ \iota_\xi(\sigma)\leq^\rho_k \iota_\xi(\tau)$
\end{itemize}
whenever $1\leq \sigma\leq\tau\leq\lambda$.
Let $J=\varphi_\xi^{-1}[[0,\lambda]]$.
\begin{enumerate}
\item
Assume $Z\subseteq J$ is finite, $\rho$-closed and  $\kappa_\alpha$-closed with $\nu_\xi\in Z$. 
Also, suppose $g:\varphi_\xi[Z]\rightarrow \lambda$ is an $\alpha$-covering of $\varphi_\xi[Z]$ with $g(0)=0$. 
 There exists a $\rho$-covering $g':Z\rightarrow \iota_\xi(\lambda)$ which is an embedding of $Z$ into ${\cal R}^{\kappa_\alpha}$ such that
 \begin{itemize}
 \item[$(*)$]
 \hspace{.4in} $g'(\mu +\chi)=(\iota_\xi\circ g\circ\iota_\xi^{-1})(\mu)+\chi$
 \end{itemize}
whenever $\mu+\chi\in Z$, $\mu$ is divisible by $\kappa_\alpha$ and $\chi<\kappa_\alpha$.
\item
Assume $Z\subseteq [0,\lambda]$ is $\alpha$-closed and $\iota_\xi(\lambda)\in J_\xi$. 
Also, suppose $Z^*\subseteq J$ is locally incompressible, $\nu_\xi\in Z^*$, $\iota_\xi[Z]\subseteq Z^*$ and $h:Z^*\rightarrow J$ is a $\rho$-covering of $Z^*$ with $h(\nu_\xi)=\nu_\xi$.
There exists an $\alpha$-covering $h':Z\rightarrow [0,\lambda]$ such that
   \begin{itemize}
   \item[$(**)$]
   \hspace{.7in} $h'(\sigma+\gamma)=(\iota_\xi^{-1}\circ h\circ\iota_\xi)(\sigma)+\gamma$
   \end{itemize}
whenever $\sigma+\gamma\in Z$, $\sigma$ is divisible by $\alpha$ and $\gamma<\alpha$.
\item
Assume $\lambda'\in[1,\lambda]$.
$\lambda'\leq^\alpha_1 \lambda+1$ iff $\lambda+1$ is in the domain of $\iota_\xi$ and $\iota_\xi(\lambda')\leq^\rho_1 \iota_\xi(\lambda+1)$.
\end{enumerate}
\end{lem}
{\bf Proof.}
Notice that $J$ is an initial segment of $\overline{J}_\xi$. 
Let $\tilde{\delta}$ and $\tilde{\gamma}$ be ordinals such that $\lambda=\alpha\cdot \tilde{\delta}+\tilde{\gamma}$ and $\tilde{\gamma}<\alpha$.
By part 13 of Lemma 8.6 of [\ref{Ca??}], $J=[\nu_\xi,\nu_\xi+\kappa_\alpha\cdot\tilde{\delta}+\kappa_{\tilde{\gamma}+1})$.

\halfblankline

{\bf (Part 1) }
By Lemma 2.4 of [\ref{Ca??}], in order to show that $(*)$ defines an embedding of $Z$ into ${\cal R}^{\kappa_\alpha}$, it suffices to show that if $\mu\in Z$ is divisible by $\kappa_\alpha$ then $g'(\mu)$ is defined and divisible by $\kappa_\alpha$.

Assume $\mu\in Z$ is divisible by $\kappa_\alpha$.
Our assumption implies $\mu\in \overline{J}'_\xi$.
Therefore, $\mu$ is in the domain of $\iota_\xi^{-1}$.
Since $\mu$ is divisible by $\kappa_\alpha$, $\iota_\xi^{-1}(\mu)$ is divisible by $\alpha$.
Moreover, $\iota_\xi^{-1}(\mu)=\varphi_\xi(\mu)$ implying $\iota_\xi^{-1}(\mu)$ is in the domain of $g$.
Since $\iota_\xi^{-1}(\mu)$ is divisible by $\alpha$ and $g$ is an $\alpha$-covering, $g(\iota_\xi^{-1}(\mu))$ is divisible by $\alpha$.
Therefore, $g'(\mu)=\iota_\xi(g(\iota_\xi^{-1}(\mu)))$ is divisible by $\kappa_\alpha$.

Since $\rho$ divides $\kappa_\alpha$, $Z$ is $\rho$-closed and $g'$ is an embedding of $Z$ into ${\cal R}^{\kappa_\alpha}$, $g'$ is clearly an embedding of $Z$ into ${\cal R}^\rho$.

To show that the range of $g'$ is contained in $\iota_\xi(\lambda)$, assume $\mu+\chi\in Z$ where $\mu$ is divisible by $\kappa_\alpha$ and $\chi<\kappa_\alpha$.
To establish $g'(\mu+\chi)<\iota_\xi(\lambda)$, it suffices to show $\varphi_\xi(g'(\mu+\chi))<\lambda$.
Let $\gamma=index(\chi)$.
Since $\varphi_\xi(\mu+\chi)=\iota_\xi^{-1}(\mu)+\gamma$, $\iota_\xi^{-1}(\mu)+\gamma$ is in the domain of $g$.
\begin{tabbing}
\hspace{.1in} $\varphi_\xi(g'(\mu+\chi))$ \ \= = \ $\varphi_\xi(g'(\mu)+\chi))$ \hspace{.5in} \= (by $(*)$) \\
\> = \ $\iota_\xi^{-1}(g'(\mu))+\gamma$ \> (since $g'(\mu)$ is divisible by $\kappa_\alpha$) \\
\>  = \ $g(\iota_\xi^{-1}(\mu))+\gamma$ \> (by $(*)$) \\
\>  = \ $g(\iota_\xi^{-1}(\mu)+\gamma)$   \> (since $\iota_\xi^{-1}(\mu)$ is divisible by $\alpha$) \\
\> $<$ \ $\lambda$ \> (since $ran(g)\subseteq \lambda$)
\end{tabbing}

To complete the proof that $g'$ is a $\rho$-covering of $Z$, assume $\mu_i+\chi_i\in Z$  where $\mu_i$ is divisible by $\kappa_\alpha$ and $\chi_i<\kappa_\alpha$ for $i=1,2$ and 
$$\mu_1+\chi_1<^\rho_k \mu_2+\chi_2$$
We will show $g'(\mu_1+\chi_1)\leq^\rho_k g'(\mu_2+\chi_2)$. Since $Z$ is $\kappa_\alpha$-closed, $\mu_i\in Z$ for $i=1,2$.

Since $g'$ is an embedding of $Z$ into ${\cal R}^{\kappa_\alpha}$,  $g'(\mu_i)$ is divisible by $\kappa_\alpha$ for $i=1,2$.

{\bf Case 1 of Part 1.} Assume $\mu_1+\chi_1=\nu_\xi$.

Since $g(0)=0$, one easily checks $g'(\nu_\xi)=\nu_\xi$.
So, we need to show $\nu_\xi\leq^\rho_k g'(\mu_2+\chi_2)$.
When $k=1$ this follows from the fact that $\nu_\xi\leq^\rho_1 \beta$ for all $\beta\in \overline{J}_\xi$.
Suppose $k=2$. 
Since $\nu_\xi\leq^\rho_2 \mu_2+\chi_2$, RTSI implies $\chi_2=0$.
Since $g'$ is an embedding of $Z$ into ${\cal R}^{\kappa_\alpha}$, $g'(\mu_2)$ is divisible by $\kappa_\alpha$.
Since $g'(\mu_2)<\iota_\xi(\lambda)$, $g'(\mu_2)\in J_\xi$.
This implies that $\nu_\xi\leq^\rho_2 \mu_2$.

{\bf Case 2 of Part 1.} Assume $\mu_1+\chi_1\not=\nu_\xi$, $\chi_1=0$ and $\chi_2\not=0$.

By RTSI, we must have $k=1$.
Let $\gamma=index(\chi_2)$.
Since $\mu_1\leq^\rho_1 \mu_2+\chi_2$ and $\kappa_\gamma\leq \chi_2$,  
$$\mu_1\leq^\rho_1 \mu_2+\kappa_\gamma$$
Since $\nu_\xi<\mu_1<\mu_2+\chi_2$, the assumptions of the lemma imply
$$\iota_\xi^{-1}(\mu_1)\leq^\alpha_1 \iota_\xi^{-1}(\mu_2+\kappa_\gamma)=\iota_\xi^{-1}(\mu_2)+\gamma$$
and $0<\iota_\xi^{-1}(\mu_1)$.
Since $\iota_\xi^{-1}(\mu_1)=\varphi_\xi(\mu_1)\in\varphi_\xi[Z]$, $\iota_\xi^{-1}(\mu_2)+\gamma=\varphi_\xi(\mu_2+\chi_2)\in\varphi_\xi[Z]$, $g$ is an $\alpha$-covering of $\varphi_\xi[Z]$ and $\alpha$ divides $\iota_\xi^{-1}(\mu_2)$, 
$$g(\iota_\xi^{-1}(\mu_1))\leq^\alpha_1 g(\iota_\xi^{-1}(\mu_2)+\gamma)=g(\iota_\xi^{-1}(\mu_2))+\gamma$$
Since $0<\iota_\xi^{-1}(\mu_1)$, $0=g(0)<g(\iota_\xi^{-1}(\mu_1))$.
The assumptions of the lemma imply
$$\iota_\xi(g(\iota_\xi^{-1}(\mu_1)))\leq^\alpha_1 \iota_\xi(g(\iota_\xi^{-1}(\mu_2))+\gamma)=\iota_\xi(g(\iota_\xi^{-1}(\mu_2)))+\kappa_\gamma$$
where the equality follows from the fact that $g(\iota_\xi^{-1}(\mu_2))$ is divisible by $\alpha$.
By  the definition of $g'$,  
$$g'(\mu_1)\leq^\rho_1g'(\mu_2)+\kappa_\gamma $$
Since $index(\chi_2)=\gamma$, $\kappa_\gamma\leq^\rho_1 \chi_2$.
 RTSI impies 
 $$g'(\mu_2)+\kappa_\gamma\leq^\rho_1 g'(\mu_2)+\chi_2$$
Therefore, 
$$g'(\mu_1)\leq^\rho_1 g'(\mu_2)+\chi_2=g'(\mu_2+\chi_2)$$

{\bf Case 3 of Part 1.} Assume $\mu_1+\chi_1\not= \nu_\xi$ and $\chi_1=\chi_2=0$.

The assumptions of the lemma and the assumption that $g$ is an $\alpha$-covering easily imply that $g'(\mu_1)\leq^\rho_k g'(\mu_2)$.

{\bf Case 4 of Part 1.} Assume $\chi_1\not=0$.

By RTSI, $\mu_1=\mu_2$ and $\chi_1\leq^\rho_k\chi_2$.
 RTSI also implies $g'(\mu_1)+\chi_1\leq^\rho_k g'(\mu_1)+\chi_2$ i.e. $g'(\mu_1+\chi_1)\leq^\rho_k g'(\mu_2+\chi_2)$.

\halfblankline

{\bf (Part 2) }
Since $\iota_\xi(\lambda)\in J_\xi$,  part 13 of Lemma 8.6 of [\ref{Ca??}] implies that $J\subseteq J_\xi$.

By parts 9 and 10 of Lemma 8.6 of [\ref{Ca??}], an element $\mu$ of $J_\xi$ is divisible by $\kappa_\alpha$ iff $\nu_\xi\leq^\rho_2 \mu$. 
Since $\nu_\xi\in Z^*$ and $h$ fixes $\nu_\xi$,  $h(\mu)$ is divisible by $\kappa_\alpha$ whenever $\mu\in Z^*$ is divisible by $\kappa_\alpha$.

By Lemma 2.4 of [\ref{Ca??}], in order to show that $(**)$ defines an embedding of $Z$ into ${\cal R}^\alpha$, it suffices to show that if $\sigma\in Z$ is divisible by $\alpha$ then $h'(\sigma)$ is defined and divisible by $\alpha$.

Assume $\sigma\in Z$ is divisible by $\alpha$.
Our assumption implies $\iota_\xi(\sigma)$ is divisible by $\kappa_\alpha$.
Moreover, the assumptions of part 2 imply $\iota_\xi(\sigma)$ is in the domain of $h$.
By the remark above, $h(\iota_\xi(\sigma))$ is divisible by $\kappa_\alpha$.
This implies that $h(\iota_\xi(\sigma))\in \overline{J}_\xi'$ and $h'(\sigma)=\iota_\xi^{-1}(h(\iota_\xi(\sigma)))$ is divisible by $\alpha$.

\halfblankline

{\bf Claim for Part 2.} If $\sigma+\gamma\in Z$ where $\sigma$ is divisible by $\alpha$ and $\gamma<\alpha$ then $h'(\sigma+\gamma)\leq \varphi_\xi(h(\iota_\xi(\sigma+\gamma)))$.

\begin{tabbing}
\hspace{.1in} $\varphi_\xi(h(\iota_\xi(\sigma+\gamma)))$ \ \= = \ $\varphi_\xi(h(\iota_\xi(\sigma)+\kappa_\gamma))$ \hspace{.5in} \= (remarks before the lemma) \\
\> $\geq$ \ $\varphi_\xi(h(\iota_\xi(\sigma))+\kappa_\gamma)$ \> (part 3 of Lemma \ref{lsrtt}) \\
\>  = \ $\iota_\xi^{-1}(h(\iota_\xi(\sigma)))+\gamma$ \> ($h(\iota_\xi(\sigma))$ is divisible by $\kappa_\alpha$) \\
\>  = $h'(\sigma)+\gamma$ \> \\
\> = $h'(\sigma+\gamma)$ \>
\end{tabbing}

\halfblankline

To see that the range of $h'$ is contained in $[0,\lambda]$, assume $\sigma+\gamma\in Z$ where $\sigma$ is divisible by $\alpha$ and $\gamma<\alpha$.
By the assumptions for part 2, $\iota_\xi(\sigma+\gamma)\in Z^*$ and $h(\iota_\xi(\sigma+\gamma))\in J$.
By definition of $J$,  $\varphi_\xi(h(\iota_\xi(\sigma+\gamma)))\in [0,\lambda]$. 
This and the claim imply $h'(\sigma+\gamma)\leq \lambda$. 

To complete the proof of part 2, assume $\sigma_i+\gamma_i\in Z$ where $\sigma_i$ is divisible by $\alpha$ and $\gamma_i<\alpha$ for $i=1,2$ and 
$$\sigma_1+\gamma_1<^\alpha_k\sigma_2+\gamma_2$$
We will show $h'(\sigma_1+\gamma_1)\leq^\alpha_k h'(\sigma_2+\gamma_2)$.
Our assumption implies that $\gamma_1=0$ and $\sigma_1\not=0$.

{\bf Case 1 of Part 2.} Assume $\gamma_2\not=0$.

By RTSI, $k=1$.
By the assumptions of the lemma, 
$$\iota_\xi(\sigma_1)\leq^\rho_1 \iota_\xi(\sigma_2+\gamma_2)=\iota_\xi(\sigma_2)+\kappa_{\gamma_2}$$
Since $h$ is a $\rho$-covering, 
$$h(\iota_\xi(\sigma_1))\leq^\rho_1 h(\iota_\xi(\sigma_2)+\kappa_{\gamma_2})$$
By part 3 of Lemma \ref{lsrtt}, $h(\iota_\xi(\sigma_2))+\kappa_{\gamma_2}   \leq h(\iota_\xi(\sigma_2)+\kappa_{\gamma_2})$.
Therefore, 
$$h(\iota_\xi(\sigma_1))\leq^\rho_1 h(\iota_\xi(\sigma_2))+\kappa_{\gamma_2}$$
Since $\sigma_1\not=0$, $\nu_\xi=\iota_\xi(0)<\iota_\xi(\sigma_1)$ implying $0\leq h(\nu_\xi)<h(\iota_\xi(\sigma_1))$.
Since $h(\iota_\xi(\sigma_2))+\kappa_{\gamma_2}$ is in the domain of $\iota_\xi^{-1}$, this implies 
$$h'(\sigma_1)\leq^\alpha_1 h'(\sigma_2+\gamma_2)$$
by the assumptions of the lemma and the definition of $h'$.

{\bf Case 2 of Part 2.} Assume $\gamma_2=0$.

The assumptions the lemma and the assumption that $h$ is a $\rho$-covering easily imply the desired conclusion.

\halfblankline

{\bf (Part 3)}

$(\Rightarrow)$ 
Assume $\lambda'\leq^\alpha_1 \lambda+1$.
Recall that $\lambda=\alpha\cdot\tilde{\delta}+\tilde{\gamma}$ where $\tilde{\gamma}<\alpha$ and  $J=[\nu_\xi,\tau)$ where $\tau=\nu_\xi+\kappa_\alpha\cdot\tilde{\delta}+\kappa_{\tilde{\gamma}+1}$.
Notice that if $\lambda+1$ were in the domain of $\iota_\xi$ then $\iota_\xi(\lambda+1)=\tau$.
By Lemma 3.1 of [\ref{Ca??}], $\lambda'$ is a limit multiple of $\alpha$ implying $\iota_\xi(\lambda')$ is a limit multiple of $\kappa_\alpha$.

\halfblankline

{\bf Claim for $(\Rightarrow)$ of part 3.} Assume $X\subseteq \iota_\xi(\lambda')$ and $Y\subseteq [\iota_\xi(\lambda'),\tau)$ are finite sets which are both $\rho$-closed and $\kappa_\alpha$-closed. There exists a $\rho$-covering of $X\cup Y$ into $\iota_\xi(\lambda')$ which fixes $X$ and is an embedding of $X\cup Y$ into ${\cal R}^{\kappa_\alpha}$.

Let $X_0$ be the set whose elements are $\nu_\xi$ along with the set of all $\mu\in X$ such that $\mu\leq^\rho_2 \mu'$ for some $\mu'\in Y$. 
By Lemma 3.9 of [\ref{Ca??}], it suffices to find a $\rho$-covering of $X_0\cup Y$ into $\iota_\xi(\lambda')$ which is an embedding of $X_0\cup Y$ into ${\cal R}^{\kappa_\alpha}$, fixes $X_0$ and maps $Y$ above $X$ (clearly, the extension of such a $\rho$-covering of $X_0\cup Y$ to $X\cup Y$ which fixes the elements of $X$ will still be an embedding into ${\cal R}^{\kappa_\alpha}$). 

By part 4 of Lemma 8.6 of [\ref{Ca??}], $\nu_\xi\leq X_0$.
By RTSI, each element of $X_0$ is divisible by $\kappa_\alpha$ implying $X_0$ is both $\rho$-closed and $\kappa_\alpha$-closed.

Since  $X_0$ and $Y$ are both $\rho$-closed and $\kappa_\alpha$-closed, $X_0\cup Y$ is both $\rho$-closed and $\kappa_\alpha$-closed.
By the remarks preceding the previous lemma, $\varphi_\xi[X_0]$ and $\varphi_\xi[Y]$ are $\alpha$-closed.

If $\beta\in X\cap\overline{J}_\xi$ then, since $\beta<\iota_\xi(\lambda')$, $\varphi_\xi(\beta)<\lambda'$.
Since $\lambda'\leq^\alpha_1\lambda$, there is an $\alpha$-covering $g$ of $\varphi_\xi[X_0]\cup \varphi_\xi[ Y]$ into $\lambda'$ which fixes $\varphi_\xi[X_0]$ and maps $\varphi_\xi[Y]$ above $\varphi_\xi(\beta)$ for all $\beta\in X\cap\overline{J}_\xi$.
Since $\nu_\xi\in X_0$, $0\in \varphi_\xi[X_0]$ and $g(0)=0$.
By part 1, there is a  $\rho$-covering $g':X_0\cup Y\rightarrow \iota_\xi(\lambda)$ which is an embedding of $X_0\cup Y$ into ${\cal R}^{\kappa_\alpha}$ such that $g'(\mu+\chi)=(\iota_\xi\circ g \circ \iota_\xi^{-1})(\mu)+\chi$ whenever $\mu+\chi\in X_0\cup Y$, $\mu$ is divisible by $\kappa_\alpha$ and $\chi<\kappa_\alpha$.

To see that $g'$ fixes the elements of $X_0$, suppose $\mu+\chi\in X_0$ where $\mu$ is divisible by $\kappa_\alpha$ and $\chi<\kappa_\alpha$.
Since all elements of $X_0$ are divisible by $\kappa_\alpha$, $\chi=0$.
Since $g$ fixes elements of $\varphi_\xi[X_0]=\iota_\xi^{-1}[X_0]$, $g'$ clearly fixes $\mu$.

To see that $g'$ maps $Y$ above $X$, it is enough to show $g'$ maps $Y$ above $X\cap \overline{J}_\xi$. Suppose $\beta\in X\cap\overline{J}_\xi$ and $\mu+\chi\in Y$ where $\mu$ is divisible by $\kappa_\alpha$ and $\chi<\kappa_\alpha$.
To show $\beta< g'(\mu+\chi)$, it suffices to show $\beta<g'(\mu)$.
Since $g'(\mu)=\iota_\xi(g(\iota_\xi^{-1}(\mu)))$, it suffices to have $\varphi_\xi(\beta)<g(\iota_\xi^{-1}(\mu))$.
Since $\iota_\xi^{-1}(\mu)=\varphi_\xi(\mu)$ (by one of the the observations before the lemma), this follows from the choice of $g$.

To see that $g'$ maps $Y$ into $\iota_\xi(\lambda')$, assume $\mu+\chi\in Y$ where $\mu$ is divisible by $\kappa_\alpha$ and $\chi<\kappa_\alpha$.
We will show $g'(\mu+\chi)<\iota_\xi(\lambda')$.
Since $\iota_\xi(\lambda')$ is a multiple of $\kappa_\alpha$ and $g'(\mu+\chi)=g'(\mu)+\chi$, it suffices to show $g'(\mu)<\iota_\xi(\lambda')$.
Since  $g'(\mu)=\iota_\xi(g(\iota_\xi^{-1}(\mu)))$, this follows from the fact that $g(\iota_\xi^{-1}(\mu))<\lambda'$.

This completes the proof of the claim.

\halfblankline

To show that $\iota_\xi(\lambda')\leq^\rho_1 \tau$, assume $X\subseteq \iota_\xi(\lambda')$ and $Y\subseteq [\iota_\xi(\lambda'),\tau)$ are finite $\rho$-closed sets. 
We will show there is a $\rho$-covering of $X\cup Y$ into $\iota_\xi(\lambda')$ which fixes $X$. 

Let $X^+$ be the union of $X$ with the collection of $\mu$ which are divisible by $\kappa_\alpha$ such that $X$ intersects $[\mu,\mu+\kappa_\alpha)$.
 Since $\kappa_\alpha$ is divisible by $\rho$, $X^+$ is both $\rho$-closed and $\kappa_\alpha$-closed.
Clearly, $X^+\subseteq \iota_\xi(\lambda')$.
 Similarly, let $Y^+$ be the union of $Y$ with the collection of $\mu$ which are divisible by $\kappa_\alpha$ such that $Y$ intersects $[\mu,\mu+\kappa_\alpha)$.
As with $X^+$,  $Y^+$ is both $\rho$-closed and $\kappa_\alpha$-closed.

By the claim, there is a $\rho$-covering of $X^+\cup Y^+$ into $\iota_\xi(\lambda')$ which fixes $X^+$. The restriction to $X\cup Y$ is the desired $\rho$-covering of $X\cup Y$.

We still need to show $\lambda+1$ is in the domain of $\iota_\xi$.
For this, it suffices to show $\tau\in \overline{J}_\xi$.

First consider the case when $\xi+1=\theta_2$. 
Since $\iota_\xi(\lambda')\leq^\rho_1 \tau$ and $\iota_\xi(\lambda')\in I_\alpha$, $\tau\in I_\alpha$ implying $\tau\in J_\xi=\overline{J}_\xi$.

Suppose $\xi+1<\theta_2$.
In this case, $\overline{J}_\xi=[\nu_\xi,\nu_{\xi+1}]$.
Argue by contradiction and assume that $\tau\not\in\overline{J}_\xi$.
This implies $\nu_{\xi+1}<\tau$.
We will derive a contradiction by showing this implies $\nu_\xi\leq^\rho_2 \nu_{\xi+1}$.

By part 14 of Lemma 8.6 of [\ref{Ca??}], $\nu_\xi \leq^\rho_2\uparrow \nu_{\xi+1}$.

To show $\nu_\xi\leq^\rho_2\downarrow$, assume $X\subseteq \nu_\xi$ and $Y\subseteq [\nu_\xi,\nu_{\xi+1})$ are finite $\rho$-closed sets. 
We will show there is a $\rho$-covering $h$ of $X\cup Y$ into $\nu_\xi$ which fixes $X$ such that $h(\zeta)\leq^\rho_1 \nu_\xi$ whenever $\zeta\in Y$ and $\zeta\leq^\rho_1\nu_{\xi+1}$.
As above, there exists finite $X^+\subseteq \nu_\xi$ containing $X$ which is both $\rho$-close and $\kappa_\alpha$-closed, and there exists finite $Y^+\subseteq [\nu_\xi,\nu_{\xi+1})$ containing $Y$ which is both $\rho$-closed and $\kappa_\alpha$-closed.
We may assume $\nu_\xi\in Y^+$.
Let $Y_1=Y^+\cap[\nu_\xi,\iota_\xi(\lambda'))$ and $Y_2=Y^+\cap [\iota_\xi(\lambda'),\nu_{\xi+1})$. 
Clearly, $Y_1$ and $Y_2$ are both $\rho$-closed and $\kappa_\alpha$-closed.
By the claim, there is a $\rho$-covering $h_1$ of $X^+\cup Y_1\cup Y_2  \cup \{\nu_{\xi+1}\}$ into $\iota_\xi(\lambda')$ which is an embedding into ${\cal R}^{\kappa_\alpha}$ and fixes $X^+\cup Y_1$.
Since $h_1$ is an embedding into ${\cal R}^{\kappa_\alpha}$, $h_1(\nu_{\xi+1})$ is divisible by $\kappa_\alpha$.
Moreover, $\nu_\xi=h_1(\nu_\xi)< h_1(\nu_{\xi+1})<\iota_\xi(\lambda')<\iota_\xi(\lambda)\leq\nu_{\xi+1}$.
Therefore, $\nu_\xi\leq^\rho_2 h_1(\nu_{\xi+1})$.
This implies there is a $\rho$-covering $h_2$ of $X^+\cup Y_1\cup h_1[Y_2]$ into $\nu_\xi$ which fixes $X^+$ such that $h_2(\zeta)\leq^\rho_1 \nu_\xi$ whenever $\zeta\in Y_1\cup h_1[Y_2]$ and $\zeta\leq^\rho_1h_1(\nu_{\xi+1})$.
The composition of $h_1$ and $h_2$ restricted to $X\cup Y$ is the desired $\rho$-covering $h$.

We have established that $\nu_\xi\leq^\rho_2\nu_{\xi+1}$ -- contradiction.

\halfblankline

$(\Leftarrow)$
Assume $\lambda+1$ is in the domain of $\iota_\xi$ and $\iota_\xi(\lambda')\leq^\rho_1 \iota_\xi(\lambda+1)$. 

To show $\lambda'\leq^\alpha_1\lambda+1$, assume $X\subseteq \lambda'$ and $Y\subseteq [\lambda',\lambda]$ are finite $\alpha$-closed sets. 
We will show there is an $\alpha$-covering of $X\cup Y$ into $\lambda'$ which fixes $X$.

Let $X_0$ consist of $0$ along with those $\sigma$ in $X$ such that $\sigma\leq^\alpha_2 \sigma'$ for some $\sigma'\in Y$.
By Lemma 3.9 of [\ref{Ca??}], it suffices to find a covering of $X_0\cup Y$ into $\lambda'$ which fixes $X_0$ and maps the elements of $Y$ above $X$.

By Lemma 3.1 of [\ref{Ca??}], each element of $X_0$ is a multiple of $\alpha$.
Therefore, each element of $\iota_\xi[X_0]$ is a multiple of $\kappa_\alpha$.
This implies that $\iota_\xi[X_0]$ is locally incompressible.
Since $0\in X_0$, $\nu_\xi=\iota_\xi(0)\in \iota_\xi[X_0]$.

By part 2 of Lemma \ref{lsrtt}, there is a finite subset $Y^*$ of $[\iota_\xi(\lambda'),\iota_\xi(\lambda+1))$ which is locally incompressible such that $\iota_\xi[Y]\subseteq Y^*$.
Since both $\iota_\xi[X_0]$ and $Y^*$ are locally incompressible, $\iota_\xi[X_0]\cup Y^*$ is locally incompressible.

Since $\iota_\xi(\lambda')\leq^\rho_1\iota_\xi(\lambda+1)$, there is a $\rho$-covering $h$ of $\iota_\xi[X_0]\cup Y^*$ into $\iota_\xi(\lambda')$ which fixes $\iota_\xi[X_0]$ and maps $Y^*$ above $\iota_\xi(\zeta)$ for each $\zeta\in X$.
By part 2, there is an $\alpha$-covering $h'$ of $X_0\cup Y$ into $[0,\lambda]$ such that $h'(\sigma+\gamma)=(\iota_\xi^{-1}\circ h \circ \iota_\xi)(\sigma)+\gamma$ whenever $\sigma+\gamma\in X_0\cup Y$, $\sigma$ is divisible by $\alpha$ and $\gamma<\alpha$.

To see that $h'$ fixes $X_0$, assume $\sigma+\gamma\in X_0$ where $\sigma$ is divisible by $\alpha$ and $\gamma<\alpha$.
Since each element of $X_0$ is a multiple of $\alpha$, $\gamma=0$.
Since $\iota_\xi(\sigma)\in \iota_\xi[X_0]$ and $\iota_\xi[X_0]$ is fixed by $h$, $\sigma$ is fixed by $h'$.

To see that $h'$ maps  $Y$ above $X$, assume $\zeta\in X$ and $\sigma+\gamma\in Y$ where $\sigma$ is divisible by $\alpha$ and $\gamma<\alpha$.
Since $\sigma\in Y$, it suffices to show $\zeta< h'(\sigma)$.
Since $\iota_\xi(\sigma)\in Y^*$, $\iota_\xi(\zeta)<h(\iota_\xi(\sigma))$ by choice of $h$.
Applying $\iota_\xi^{-1}$, $\zeta< \iota_\xi^{-1}(h(\iota_\xi(\sigma)))=h'(\sigma)$.

To see that $h'$ maps $Y$ into $\lambda'$, assume $\sigma+\gamma\in Y$ where $\sigma$ is divisible by $\alpha$ and $\gamma<\alpha$.
Since $\lambda'$ is divisible by $\kappa_\alpha$ and $h'(\sigma+\gamma)=h'(\sigma)+\gamma$, to show $h'(\sigma+\gamma) <\lambda'$ it suffices to establish that $h'(\sigma)<\lambda'$. 
This follows from the definition of $h'$ and the fact that $h(\iota_\xi(\sigma))<\iota_\xi(\lambda')$. 
\qed


\begin{thm}
{\bf (Order Reduction Theorem)}
Assume $\xi<\theta_2$.
If $1\leq \lambda' < \lambda \in dom(\iota_\xi)$ then
\begin{itemize}
\item[$(\star)$]
\hspace{1.1in} $\lambda'\leq^\alpha_k \lambda \ \ \ \ \Longleftrightarrow \ \ \ \iota_\xi(\lambda')\leq^\rho_k \iota_\xi(\lambda)$
\end{itemize}
for $k=1,2$.
\end{thm}
{\bf Proof.}
Let $K$ be the collection of ordinals $\lambda$ in the domain of $\iota_\xi$ such that for $k=1,2$ and all $\lambda'\in [1,\lambda)$, $(\star)$ holds.
We will show that $K=dom(\iota_\xi)$ by induction.

Assume $\lambda\in dom(\iota_\xi)$ and $\lambda'\in K$ for all $\lambda'<\lambda$.

First, suppose $k=1$ and assume $\lambda'\in [1,\lambda)$.
If $\lambda$ is a limit ordinal then $(\star)$ follows from the induction hypothesis using the continuity of $\iota_\xi$ and part 2 of Lemma 3.8 of [\ref{Ca??}]. On the other hand, if $\lambda$ is a successor ordinal then $(\star)$ follows from part 3 of Lemma \ref{lorA}.

We now consider the case $k=2$. 
Assume $\lambda'\in[1,\lambda)$. 

\halfblankline

($\Rightarrow$)
Assume $\lambda'\leq^\alpha_2\lambda$.

By Lemma 3.1 of [\ref{Ca??}] and part 1 of Lemma 3.6 of [\ref{Ca??}], both $\lambda'$ and $\lambda$ are limit multiples of $\alpha$.
Therefore, both $\iota_\xi(\lambda')$ and $\iota_\xi(\lambda)$ are limit multiples of $\kappa_\alpha$.

\halfblankline

We first show that $\iota_\xi(\lambda')\leq^\rho_2\downarrow \iota_\xi(\lambda)$.
The proof is similar to the forward direction of part 3 of Lemma \ref{lorA}.

Assume $X\subseteq \iota_\xi(\lambda')$ and $Y\subseteq [\iota_\xi(\lambda'),\iota_\xi(\lambda))$ are finite $\rho$-closed sets. 
Let $X_0$ be the set whose elements are $\nu_\xi$ along with all $\mu\in X$ such that $\mu\leq^\rho_2\mu'$ for some $\mu'\in Y$.
By Lemma 3.9 of [\ref{Ca??}], it suffices to find a $\rho$-covering $f$ of $X_0\cup Y$ into $\iota_\xi(\lambda')$ which fixes $X_0$ and maps $Y$ above $X$ such that $f(\mu)\leq^\rho_1 \iota_\xi(\lambda')$ whenever $\mu\in Y$ and $\mu\leq^\rho_1\iota_\xi(\lambda)$.

By part 4 of Lemma 8.6 of [\ref{Ca??}], $\nu_\xi\leq X_0$.
By RTSI, each element of $X_0$ is divisible by $\kappa_\alpha$ implying $X_0$ is both $\rho$-closed and $\kappa_\alpha$-closed.

Let $Y^+$ be the union of $Y$ with the collection of $\mu$ which are divisible by $\kappa_\alpha$ such that $Y$ intersects $[\mu,\mu+\kappa_\alpha)$.
 Since $\kappa_\alpha$ is divisible by $\rho$, $Y^+$ is both $\rho$-closed and $\kappa_\alpha$-closed.
Since $\iota_\xi(\lambda')$ is divisible by $\kappa_\alpha$, $Y^+\subseteq [\iota_\xi(\lambda'),\iota_\xi(\lambda))$.

Without loss of generality, we may replace $Y$ by $Y^+$

 Since  $\lambda'\leq^\alpha_2\lambda$, there is a covering $g$ of $\varphi_\xi[X_0]\cup \varphi_\xi[Y^+]$ into $\lambda'$ such that $g$ fixes $\varphi_\xi[X_0]$, $\varphi_\xi(\beta)<g[\varphi_\xi[Y^+]]$  for each $\beta\in X\cap\overline{J}_\xi$ and $g(\sigma)\leq^\alpha_1 \lambda'$ whenever $\sigma\in \varphi_\xi[Y]$ and $\sigma\leq^\alpha_1\lambda$.
 
Since $\nu_\xi\in X_0$, $0=\varphi_\xi(\nu_\xi)\in \varphi_\xi[X_0]$ and  $g(0)=0$.
 
By part 1 of Lemma \ref{lorA},  there is a $\rho$-covering $g'$ of $X_0\cup Y^+$ such that $g'(\mu+\chi)= (\iota_\xi^{-1}\circ g \circ \iota_xi)(\mu)+\chi$ whenever $\mu+\chi\in X_0\cup Y^+$, $\mu$ is divisible by $\kappa_\alpha$ and $\chi<\kappa_\alpha$.

By a straightforward argument similar to that in the forward direction of part 3 of Lemma \ref{lorA}, $g'$ maps $X_0\cup Y^+$ into $\iota_\xi(\lambda')$, fixes $X$ and maps $Y$ above $X$.

Suppose $\mu\in Y^+$ and $\mu\leq^\rho_1\iota_\xi(\lambda)$.
Since $\iota_\xi(\lambda)$ is a multiple of $\kappa_\alpha$, $\mu$ must be a multiple of $\kappa_\alpha$ by RTSI. 
Therefore, $\mu$ is in the range of $\iota_\xi$ and, by the case $k=1$, $\iota_\xi^{-1}(\mu)\leq^\alpha_1 \lambda$.
Therefore, $g(\iota_\xi^{-1}(\mu))\leq^\alpha_1\lambda'$ implying $g'(\mu)\leq^\rho_1 \iota_\xi(\lambda')$ by the induction hypothesis.

The restriction of $g'$ to $X_0\cup Y$ is the desired $\rho$-covering $f$.

\halfblankline

We now show $\iota_\xi(\lambda')\leq^\rho_2\uparrow \iota_\xi(\lambda)$.

Assume $X$ is a finite $\rho$-closed subset of $\iota_\xi(\lambda')$ and $\cal Y$ is a collection of finite $\rho$-closed sets which is cofinal in $\iota_\xi(\lambda')$ such that $X<Y$ for $Y\in \cal Y$ and $X\cup Y_1\cong_\rho X\cup Y_2$ for all $Y_1,Y_2\in \cal Y$.
Also, suppose $\beta<\iota_\xi(\lambda)$.
We will show there is $Y\in \cal Y$ and a covering of $X\cup Y$ into $\iota_\xi(\lambda)$ which fixes $X$ and maps $Y$ above $\beta$. 
This suffices by Lemma 3.4 of [\ref{Ca??}].

Since $\lambda$ is a limit multiple of $\alpha$ and $\iota_\xi$ is continuous, we may assume there is $\beta'<\lambda$ which is a multiple of $\alpha$ such that $\beta=\iota_\xi(\beta')$. 
Moreover, we may assume $\lambda'<\beta'$.
Our assumptions imply  $\beta$ is a multiple of $\kappa_\alpha$.

Let $X_0$ be the set whose elements are $\nu_\xi$ along with those $\mu\in X$ such that for each $Y\in \cal Y$ there exists $\mu'\in Y$ such that $\mu\leq^\rho_2 \mu'$.
As before, $\nu_\xi\leq X_0$. Also, each element of $X_0$ is divisible by $\kappa_\alpha$ implying $X_0$ is locally incompressible.
In particular, $X_0$ is both $\rho$-closed and $\kappa_\alpha$-closed.

Since $\iota_\xi(\lambda')$ is a limit multiple of $\kappa_\alpha$, there is $\nu<\iota_\xi(\lambda')$ which is a multiple of $\kappa_\alpha$ such that $X_0<\nu$. 
We may assume $\nu\leq Y$ for $Y\in \cal Y$.

By Lemma 3.10 of [\ref{Ca??}], it suffices to show that there exists $Y\in \cal Y$ and a covering of $X_0\cup Y$ into $\iota_\xi(\lambda)$ which fixes $X_0$ and maps $Y$ above $\beta$. 

For $Y\in \cal Y$, let $M_Y$ be the set of $\mu$ such that $\mu$ is divisible by $\kappa_\alpha$ and $Y$ intersects $[\mu,\mu+\kappa_\alpha)$.

Since $X\cup Y_1\cong_\rho X\cup Y_2$ for $Y_1,Y_2\in \cal Y$, $rem^\rho[Y]$ and $card(Y)$ are constant as $Y$ varies over $\cal Y$. 
By Lemma \ref{lorA0} there is a $\rho$-incompressible subset $R$ of $\kappa_\alpha$ and ${\cal Y}'\subseteq {\cal Y}$ such that  $0\in R$, ${\cal Y}'$ is cofinal in $\iota_\xi(\lambda')$,    and there is a $\rho$-covering of $X_0\cup Y$ into $X_0\cup (M_Y+R)$ for all $Y\in {\cal Y}'$.

Since $0\in R$, $M_Y+R$ is locally incompressible for $Y\in {\cal Y}'$. In particular, $M_Y+R$ is both $\rho$-closed and $\kappa_\alpha$-closed for $Y\in {\cal Y}'$.

Assume $Y\in {\cal Y}'$.
Since $X_0$ and $M_Y+R$ are $\kappa_\alpha$-closed, $\varphi_\xi[X_0]\cup \varphi_\xi[M_Y+R]$ is $\alpha$-closed.
Moreover, $\varphi_\xi[M_Y+R]=\varphi_\xi[M_Y]+index[R]=\iota_\xi^{-1}[M_Y]+index[R]$ implying $rem^\alpha[\varphi_\xi[M_Y+R]]=index[R]$. 

Since $rem^\alpha[\varphi_\xi[M_Y+R]]$ is fixed as $Y$ ranges over ${\cal Y}'$, Lemma 3.3 of [\ref{Ca??}] implies there is a subcollection ${\cal Y}''$ of ${\cal Y}'$ which is cofinal in $\iota_\xi(\lambda')$ such that $\varphi_\xi[X_0]\cup \varphi_\xi[M_{Y_1}+R]\cong_\alpha \varphi_\xi[X_0]\cup \varphi_\xi[M_{Y_2}+R]$ for all $Y_1,Y_2\in {\cal Y}''$.

Since $\iota_\xi(\lambda')$ is a limit multiple of $\kappa_\alpha$ and ${\cal Y}''$ is cofinal in $\iota_\xi(\lambda')$, we clearly have that $M_Y$ $(Y\in {\cal Y}'')$ is cofinal in $\iota_\xi(\lambda')$.
This implies that $\varphi_\xi[M_Y+R]$ $(Y\in {\cal Y}'')$ is cofinal in $\lambda'$.

Assume $Y\in {\cal Y}''$. 
Since $\lambda'\leq^\alpha_2\lambda$, there is an $\alpha$-covering $g$ of $\varphi_\xi[X_0]\cup \varphi_\xi[M_Y+R]$ into $\lambda$ which fixes $\varphi_\xi[X_0]$ such that $\beta'<g[\varphi_\xi[M_Y+R]]$. 
Since $\nu_\xi\in X_0$, $0=\varphi_\xi(\nu_\xi)\in \varphi_\xi[X_0]$. 
Therefore, $g(0)=0$.
By part 1 of Lemma \ref{lorA}, there is a $\rho$-covering $g'$  of $X_0\cup (M_Y+R)$ into $\iota_\xi(\lambda)$ such that $g'(\mu+\chi)=(\iota_\xi \circ g \circ \iota_\xi^{-1})(\mu)+\chi$ whenever $\mu+\chi\in X_0\cup (M_Y+R)$, $\mu$ is divisible by $\kappa_\alpha$ and $\chi<\kappa_\alpha$.

A straightforward argument shows $g'$ fixes each element of $X_0$, being divisible by $\kappa_\alpha$, and $g'$ maps each element of $M_Y$, hence of $M_Y+R$, above $\iota_\xi(\beta')$.


\halfblankline

$(\Leftarrow)$
Assume $\iota_\xi(\lambda')\leq^\rho_2\iota_\xi(\lambda)$. 
By the case $k=1$, $\lambda'\leq^\alpha_1 \lambda$.
By Lemma 3.1 of [\ref{Ca??}], $\lambda'$ is a limit multiple of $\alpha$ implying $\iota_\xi(\lambda')$ is a limit multiple of $\kappa_\alpha$.
By RTSI, $\iota_\xi(\lambda)$ is a multiple of $\kappa_\alpha$ implying $\lambda$ is a multiple of $\alpha$.

\halfblankline

We first show $\lambda'\leq^\alpha_2\downarrow \lambda$.
The proof is similar to that of the reverse direction of part 3 of Lemma \ref{lorA}.

Assume $X\subseteq \lambda'$ and $Y\subseteq [\lambda',\lambda)$ are finite $\alpha$-closed sets. 
We will show there is an $\alpha$-covering $h'$ of $X\cup Y$ in $\lambda'$ which fixes $X$ such that $h'(\sigma)\leq^\alpha_1 \lambda'$ whenever $\sigma \in Y$ and $\sigma\leq^\alpha_1\lambda$. 

Let $X_0$ be the set whose elements are $0$ along with those $\sigma$ in $X$ such that $\sigma\leq^\alpha_2 \sigma'$ for some $\sigma'\in Y$.
By Lemma 3.9 of [\ref{Ca??}], it suffices to find a covering $h'$ of $X_0\cup Y$ into $\lambda'$ which fixes $X_0$ and maps the elements of $Y$ above $X$ such that $h'(\sigma)\leq^\alpha_1\lambda'$ whenever $\sigma\in Y$ and $\sigma\leq^\alpha_1 \lambda$.

By Lemma 3.1 of [\ref{Ca??}], each element of $X_0$ is a multiple of $\alpha$.
Therefore, each element of $\iota_\xi[X_0]$ is a multiple of $\kappa_\alpha$.
This implies that $\iota_\xi[X_0]$ is locally incompressible.
Since $0\in X_0$, $\nu_\xi=\iota_\xi(0)\in \iota_\xi[X_0]$.

By part 2 of Lemma \ref{lsrtt}, there is a finite subset $Y^*$ of $[\iota_\xi(\lambda'),\iota_\xi(\lambda))$ which is locally incompressible such that $\iota_\xi[Y]\subseteq Y^*$.
Since both $\iota_\xi[X_0]$ and $Y^*$ are locally incompressible, $\iota_\xi[X_0]\cup Y^*$ is locally incompressible.

Since $\iota_\xi(\lambda')\leq^\rho_2\iota_\xi(\lambda)$, there is a $\rho$-covering $h$ of $\iota_\xi[X_0]\cup Y^*$ into $\iota_\xi(\lambda')$ which fixes $\iota_\xi[X_0]$ and maps $Y^*$ above $\iota_\xi(\zeta)$ for each $\zeta\in X$ such that $h(\mu)\leq^\rho_1\iota_\xi(\lambda')$ whenever $\mu\in Y^*$ and $\mu\leq^\rho_1 \iota_\xi(\lambda)$.
By part 2 of Lemma \ref{lorA}, there is an $\alpha$-covering $h'$ of $X_0\cup Y$ into $\lambda$ such that $h'(\sigma+\gamma)=(\iota_\xi^{-1}\circ h \circ \iota_\xi)(\sigma)+\gamma$ whenever $\sigma+\gamma\in X_0\cup Y$, $\sigma$ is divisible by $\alpha$ and $\gamma<\alpha$.

As in the proof of the reverse direction of part 3 of Lemma \ref{lorA}, $h'$ is a covering of $X_0\cup Y$ into $\lambda'$ which fixes $X_0$ and maps $Y$ above $X$. 

Assume $\sigma\in Y$ and $\sigma\leq^\alpha_1\lambda$.
By the case $k=1$, $ \iota_\xi(\sigma)\leq^\rho_1 \iota_\xi(\lambda)$ implying $h(\iota_\xi(\sigma))\leq^\rho_1\iota_\xi(\lambda')$.
By the induction hypothesis and the definition of $h'$, $h'(\sigma)\leq^\alpha_1 \lambda'$.

\halfblankline

We now show that $\lambda'\leq^\alpha_2\uparrow \lambda$.

Assume $X$ is a finite $\alpha$-closed subset of $\lambda'$ and $\cal Y$ is a collection of nonempty finite $\alpha$-closed subsets of $\lambda'$ which is cofinal in $\lambda'$ such that $X<Y$ for  all $Y\in \cal Y$ and $X\cup Y_1 \cong_\alpha X\cup Y_2$ whenever $Y_1,Y_2\in \cal Y$.
Also, suppose that $\beta<\lambda$.
We will show that there is $Y\in \cal Y$ and an $\alpha$-covering  of $X\cup Y$ into $\lambda$ which fixes $X$ and maps $Y$ above $\beta$.
We may assume that $\lambda'\leq \beta$.

Let $X_0$ be the set whose elements are $0$ along with those $\sigma\in X$ such that for each $Y\in \cal Y$ there exists $\sigma'\in Y$ such that $\sigma\leq^\alpha_\rho\sigma'$.
By Lemma 3.1 of [\ref{Ca??}], each element of $X_0$ is a multiple of $\alpha$.
Therefore, each element of $\iota_\xi[X_0]$ is divisible by $\kappa_\alpha$ implying that $\iota_\xi[X_0]$ is locally incompressible. 
Also, $\nu_\xi=\iota_\xi(0)\in \iota_\xi[X_0]$.

By Lemma 3.10 of [\ref{Ca??}], it suffices to show that there exists $Y\in \cal Y$ and a covering of $X_0\cup Y$ into $\iota_\xi(\lambda)$ which fixes $X_0$ and maps $Y$ above $\beta$. 

Since $X\cup Y_1\cong_\alpha X\cup Y_2$ for all $Y_1,Y_2\in {\cal Y}$, $rem^\alpha[Y_1]=rem^\alpha[Y_2]$ for all $Y_1,Y_2\in \cal Y$. 
Let $R$ be the common value of $rem^\alpha[Y]$ as $Y$ ranges of $\cal Y$.
Since each element of $\cal Y$ is $\alpha$-closed, $0\in R$.

Assume $Y\in \cal Y$. 
Let $M_Y$ consist of the elements of $Y$ which are divisible by $\alpha$.
Since $Y$ is $\alpha$-closed, $Y \subseteq M_Y+R$.
Since $0\in R$, $M_Y+R$ is $\alpha$-closed.

Notice that $\iota_\xi[M_Y+R]= \iota_\xi[M_Y]+\{\kappa_\gamma \, | \, \gamma \in R\}$ and the elements of $\iota_\xi[M_Y]$ are divisible by $\kappa_\alpha$ for $Y\in \cal Y$.
Since $0\in R$, $0=\kappa_0\in \{\kappa_\gamma \, | \, \gamma \in R\}$.

By Lemma 6.7 of [\ref{Ca??}], there is a $\rho$-incompressible set $R^*$ which contains $\{\kappa_\gamma \, | \, \gamma\in R\}$ such that $index(max(R^*))$ is the largest element of $R$.
Since $0\in \{\kappa_\gamma \, | \, \gamma \in R\}$, $0\in R^*$.

Assume $Y\in \cal Y$.
By choice of $R^*$, $\iota_\xi[M_Y+R]\subseteq \iota_\xi[M_Y]+R^*$. 
Since the elements of $\iota_\xi[M_Y]$ are divisible by $\kappa_\alpha$ and $0\in R^*$, $\iota_\xi[M_Y]+R^*$ is locally incompressible.
Since both $\iota_\xi[X_0]$ and $\iota_\xi[M_Y]+R^*$ are locally incompressible, $\iota_\xi[X_0]\cup (\iota_\xi[M_Y]+R^*)$ is locally incompressible.
Since $\iota_\xi(\lambda')$ is divisible by $\kappa_\alpha$ and $M_Y<\lambda'$, $\iota_\xi[M_Y]+R^*< \iota_\xi(\lambda')$.
Using the fact that the elements of  $\iota_\xi[M_Y]$ are divisible by $\kappa_\alpha$ again, $rem^\rho[\iota_\xi[M_Y]+R^*] =rem^\rho[R^*]$.

Since $rem^\rho[\iota_\xi[M_Y]+R^*]$ is fixed as $Y$ varies over $\cal Y$,  Lemma 3.3 of [\ref{Ca??}], implies there is a subcollection ${\cal Y}'$ of $\cal Y$ which is cofinal in $\lambda'$ such that $\iota_\xi[X_0] \cup (\iota_\xi[M_{Y_1}]+ R^*) \cong_\rho \iota_\xi[X_0]\cup (\iota_\xi[M_{Y_2}]+ R^*)$ for all $Y_1,Y_2\in {\cal Y}'$.
Since $min(M_Y)\in Y$ for $Y\in {\cal Y}$, $min(M_Y)$ $(Y\in{\cal Y}')$ is cofinal in $\lambda'$.
By the continuity of $\iota_\xi$, $\iota_\xi(min(M_Y))$ $(Y\in{\cal Y}')$ is cofinal in $\iota_\xi(\lambda')$.
Therefore, $\iota_\xi[M_Y]+R^*$ $(Y\in{\cal Y}')$ is cofinal in $\iota_\xi(\lambda')$.

Fix $Y\in {\cal Y}'$. Since $\iota_\xi(\lambda')\leq^\rho_2 \iota_\xi(\lambda)$, there is a $\rho$-covering of $\iota_\xi[X_0] \cup (\iota_\xi[M_Y]+ R^*)$ into $\iota_\xi(\lambda)$ which fixes $\iota_\xi[X_0]$ and maps $\iota_\xi[M_Y]+ R^*$ above $\iota_\xi(\beta)$.
Since $\nu_\xi\in \iota_\xi[X_0]$, $h(\nu_\xi)=\nu_\xi$. 
By part 2 of Lemma \ref{lorA}, there is an $\alpha$-covering $h'$ of $X_0\cup (M_Y+R)$ into $\lambda$ such that $h'(\sigma+\gamma)= (\iota_\xi^{-1}\circ h \circ \iota_\xi)(\sigma)+\gamma$ whenever $\sigma+\gamma\in X_0\cup (M_Y+R)$, $\sigma$ is divisible by $\alpha$ and $\gamma<\alpha$.

Since the elements of $X_0$ are divisible by $\alpha$, an easy argument shows $X_0$ is fixed by $h'$.

To show that $h'$ maps $M_Y+R$ above $\beta$, it suffices to show $h'$ maps the elements of $M_Y$ above $\beta$. 
This is straightforward from the definition of $h'$ and the fact that the elements of $M_Y$ are divisible by $\alpha$ and $h$ maps $\iota_\xi[M_Y]$ above $\iota_\xi(\beta)$.
\qed

\section{The Second Recurrence Theorem for \leqtwo}
\label{ssecondrtwo}

Recall that we are assuming $\rho$ is an arbitrary additively indecomposable ordinal and $\alpha$ is an epsilon number greater than $\rho$.


\begin{lem} 
\label{lorB}
If $\eta_1<\eta_2$ and $\omega^{\eta_2}<\theta_2$ then $J_{\omega^{\eta_1}}$ is $\rho$-isomorphic to a proper initial segment of $J_{\omega^{\eta_2}}$. 
\end{lem}
{\bf Proof.}
Define $\lambda_i$ so that $J_{\omega^{\eta_i}}=[\nu_{\omega^{\eta_i}},\nu_{\omega^{\eta_i}}+\lambda_i)$ for $i=1,2$.
Notice that $\lambda_2$ may be $\infty$.
Let $\lambda$ be the minimum of $\lambda_1$ and $\lambda_2$.

\halfblankline

{\bf Claim.} 
The map from $[\nu_{\omega^{\eta_1}},\nu_{\omega^{\eta_1}}+\lambda)$ to $[\nu_{\omega^{\eta_2}},\nu_{\omega^{\eta_2}}+\lambda)$ given by
$$\nu_{\omega^{\eta_1}}+\beta \mapsto \nu_{\omega^{\eta_2}}+\beta$$
is a $\rho$-isomorphism.

As substructures of \calRrho, both $[\nu_{\omega^{\eta_1}},\nu_{\omega^{\eta_1}}+\lambda)$ and $[\nu_{\omega^{\eta_2}},\nu_{\omega^{\eta_2}}+\lambda)$ are isomorphic to initial segments of \calRrho\ implying $[\nu_{\omega^{\eta_1}},\nu_{\omega^{\eta_1}}+\lambda)$ is isomorphic to $[\nu_{\omega^{\eta_2}},\nu_{\omega^{\eta_2}}+\lambda)$ as substructures of \calRrho. 

Suppose $\kappa_\alpha\cdot\delta_j+\chi_j<\lambda$ and $\chi_j<\kappa_\alpha$ for $j=1,2$.
We need to show 
 $$\nu_{\omega^{\eta_1}}+\kappa_\alpha\cdot\delta_1+\chi_1 \leq^\rho_k \nu_{\omega^{\eta_1}}+\kappa_\alpha\cdot\delta_2+\chi_2$$
 iff
 $$\nu_{\omega^{\eta_2}}+\kappa_\alpha\cdot\delta_1+\chi_1 \leq^\rho_k \nu_{\omega^{\eta_2}}+\kappa_\alpha\cdot\delta_2+\chi_2$$
 for $k=1,2$.
This follows from RTSI when $\chi_1\not=0$ and from the Order Reduction Theorem when $\chi_1=\chi_2=0$.
Assume $\chi_1=0$ and $\chi_2\not=0$.
By RTSI, both sides of the equivalence fail if $k=2$.
So, we may assume $k=1$.
By RTSI again, $\nu_{\omega^{\eta_i}}+\kappa_\alpha\cdot\delta_2 +\kappa_\gamma\leq^\rho_1 \nu_{\omega^{\eta_i}}+\kappa_\alpha\cdot\delta_2+\chi_2$ for $i=1,2$ where $\gamma=index(\chi_2)$.
Therefore, $\nu_{\omega^{\eta_i}}+\kappa_\alpha\cdot\delta_1\leq^\rho_1 \nu_{\omega^{\eta_i}}+\kappa_\alpha\cdot\delta_2 + \chi_2$ is equivalent to $\nu_{\omega^{\eta_i}}+\kappa_\alpha\cdot\delta_1\leq^\rho_1 \nu_{\omega^{\eta_i}}+\kappa_\alpha\cdot\delta_2 + \kappa_\gamma$ for $i=1,2$.
The Order Reduction Theorem implies that $\nu_{\omega^{\eta_i}}+\kappa_\alpha\cdot\delta_1\leq^\rho_1 \nu_{\omega^{\eta_i}}+\kappa_\alpha\cdot\delta_2+\kappa_\gamma$ iff $\alpha\cdot\delta_1\leq^\alpha_1 \alpha\cdot\delta_2 +\gamma$ for $i=1,2$.
The desired equivalence follows.

\halfblankline

Argue by contradiction and assume $\lambda_2\leq \lambda_1$.

Notice that $J_{\omega^{\eta_1}}=[\nu_{\omega^{\eta_1}},\nu_{\omega^{\eta_1}+1})$ implying $\nu_{\omega^{\eta_1}}+\lambda_1=\nu_{\omega^{\eta_1}+1}$ and $\lambda_1$ is divisible by $\kappa_\alpha$.

First, assume that $\omega^{\eta_2}+1=\theta_2$.
In this case, $J_{\omega^{\eta_2}}=\overline{J}_{\omega^{\eta_2}}=[\nu_{\omega^{\eta_2}},max(I_\alpha)]$. 
We will derive a contradiction by showing $\nu_{\omega^{\eta_2}}\leq^\rho_1 max(I_\alpha)+1$.

Assume $X\subseteq \nu_{\omega^{\eta_2}}$ and $Y\subseteq J_{\omega^{\eta_2}}$ are finite $\rho$-closed sets. 
We will show there is a $\rho$-covering of $X\cup Y$ into $\nu_{\omega^{\eta_2}}$ which fixes $X$. 
By part  4 of Lemma 8.6 of [\ref{Ca??}], $X\not\leq^\rho_2 Y$.
By Lemma 3.9 of [\ref{Ca??}], it will suffice to find a $\rho$-covering of $Y$ into $\nu_{\omega^{\eta_2}}$ which maps $Y$ above $X$.
By the claim, $J_{\omega^{\eta_1}}$ contains a $\rho$-closed set which is $\rho$-isomorphic to $Y$. 
By the continuity of $\xi\mapsto \nu_\xi$, there exists $\eta<\omega^{\eta_2}$ such that $X<\nu_\eta$. 
Since $\eta,\omega^{\eta_1}<\omega^{\eta_2}$, $\eta+\omega^{\eta_1}<\omega^{\eta_2}$.
By the First Recurrence Theorem for $\leqtwo$, $J_{\eta+\omega^{\eta_1}}\cong_\rho J_{\omega_{\eta_1}}$.
Therefore, $J_{\eta+\omega^{\eta_1}}$ contains a $\rho$-closed set $\tilde{Y}$ which is $\rho$-isomorphic to $Y$.
Since $X<\nu_\eta<\tilde{Y}<\nu_{\omega^{\eta_2}}$, the $\rho$-isomorphism of $Y$ and $\tilde{Y}$ is the desired $\rho$-covering of $Y$.

Now, suppose $\omega^{\eta_2}+1<\theta_2$. 
We will derive a contradiction by showing $\nu_{\omega^{\eta_2}}\leq^\rho_2 \nu_{\omega^{\eta_2}+1}$.
In this case, $J_{\omega^{\eta_2}}=[\nu_{\omega^{\eta_2}},\nu_{\omega^{\eta_2}+1})$ and $\nu_{\omega^{\eta_2}+1}=\nu_{\omega^{\eta_2}}+ \lambda_2$.
Therefore, $\lambda_2$ is divisible by $\kappa_\alpha$.

By part 14 of Lemma 8.6 of [\ref{Ca??}], $\nu_{\omega^{\eta_2}}\leq^\rho_2\uparrow \nu_{\omega^{\eta_2}+1}$.

To show $\nu_{\omega^{\eta_2}}\leq^\rho_2\downarrow \nu_{\omega^{\eta_2}+1}$, assume $X\subseteq \nu_{\omega^{\eta_2}}$ and $Y\subseteq J_{\omega^{\eta_2}}$ are finite $\rho$-closed sets. 
We will show there is a $\rho$-covering $h$ of $X\cup Y$ into $\nu_{\omega^{\eta_2}}$ which fixes $X$ such that $h(\beta)\leq^\rho_1\nu_{\omega^{\eta_2}}$ whenever $\beta\in Y$ and $\beta\leq^\rho_1 \nu_{\omega^{\eta_2}+1}$.
By part 4 of Lemma 8.6 of [\ref{Ca??}], $X\not\leq^\rho_2 Y$.
By Lemma 3.9 of [\ref{Ca??}], it will suffice to find a $\rho$-covering $h$ of $Y$ into $\nu_{\omega^{\eta_2}}$ which maps $Y$ above $X$ such that $h(\beta)\leq^\rho_1\nu_{\omega^{\eta_2}}$ whenever $\beta\in Y$ and $\beta\leq^\rho_1 \nu_{\omega^{\eta_2}+1}$.
Let $R$ satisfy $Y=\nu_{\omega^{\eta_2}}+R$.
Clearly, $R$ is $\rho$-closed.
Choose $\eta<\omega^{\eta_2}$ such that $X<\nu_\eta$.
Let $\xi=\eta+\omega^{\eta_1}$.
Since $\eta,\omega^{\eta_1}<\omega^{\eta_2}$, $\xi<\omega^{\eta_2}$.
By the First Recurrence Theorem for $\leqtwo$, $J_{\omega^{\eta_1}}\cong_\rho J_\xi$.
This implies $\nu_{\xi+1}=\nu_\xi+\lambda_1$.
By the claim, $J_{\omega^{\eta_2}}$ is $\rho$-isomorphic to an initial segment of $J_\xi$.
Therefore, $\nu_\xi+R\cong_\rho \nu_{\omega^{\eta_2}}+R=Y$.
Also, $X<\nu_\xi+R<\nu_{\xi+1}\leq \nu_{\omega^{\eta_2}}$.

We consider two cases.

First, suppose $\lambda=\lambda_1=\lambda_2$.

We claim the $\rho$-isomorphism $\nu_{\omega^{\eta_2}}+\zeta \mapsto \nu_\xi +\zeta$ is the desired $\rho$-covering of $Y$. 
To establish this, suppose $\zeta\in R$ and $\nu_{\omega^{\eta_2}}+\zeta\leq^\rho_1 \nu_{\omega^{\eta_2}+1}$.
We will show that $\nu_\xi +\zeta \leq^\rho_1 \nu_{\omega^{\eta_2}}$.
Since $\nu_{\omega^{\eta_2}}+\zeta\leq^\rho_1 \nu_{\omega^{\eta_2}+1}=\nu_{\omega^{\eta_2}}+\lambda$, $\nu_{\omega^{\eta_2}}+\zeta\leq^\rho_1 \nu_{\omega^{\eta_2}}+\zeta'$ whenever $\zeta\leq\zeta'<\lambda_2$.
Since $J_{\omega^{\eta_2}}$ is $\rho$-isomorphic to an initial segment of $J_\xi$, $\nu_\xi+\zeta\leq^\rho_1 \nu_\xi+\zeta'$ whenever $\zeta\leq\zeta'<\lambda_2$.
By part 2 of Lemma 3.8 of [\ref{Ca??}], $\nu_\xi+\zeta\leq^\rho_1 \nu_\xi+\lambda=\nu_{\xi+1}\leq^\rho_1 \nu_{\omega^{\eta_2}}$.

Now, consider the case when $\lambda_2<\lambda_1$.
Since $\lambda_2$ is divisible by $\kappa_\alpha$, so is $\nu_\xi+\lambda_2$.
Since $\nu_\xi+\lambda_2<\nu_\xi+\lambda_1=\nu_{\xi+1}$, this implies $\nu_\xi\leq^\rho_2 \nu_\xi+\lambda_2$.
Therefore, there is a $\rho$-covering $g$ of $\nu_\xi+R$ into $\nu_\xi$ which maps $\nu_\xi+R$ above $X$ such that $g(\nu_\xi +\zeta) \leq^\rho_1 \nu_\xi$ whenever $\zeta\in R$. 
Let $h$ be the function on $Y$ defined by $h(\nu_{\omega^{\eta_2}}+\zeta)=g(\nu_\xi+\zeta)$.
Clearly, $h$ is a $\rho$-covering of $Y$.

To complete the proof of this case, assume $\zeta\in R$ and $\nu_{\omega^{\eta_2}}+\zeta\leq^\rho_1 \nu_{\omega^{\eta_2}+1}$.
By an argument similar to the previous case, $\nu_\xi+\zeta \leq^\rho_1 \nu_\xi+\lambda_2$.
Therefore, $h(\nu_{\omega^{\eta_2}}+\zeta)=g(\nu_\xi+\zeta)\leq^\rho_1\nu_\xi\leq^\rho_1 \nu_{\omega^{\eta_2}}$.
\qed


\begin{lem} 
\label{lorB1}
Assume $\omega^\eta<\theta_2$.
If $\eta'<\eta$ then $dom(\iota_{\omega^{\eta'}})\subseteq dom(\iota_{\omega^\eta})$.
\end{lem}
{\bf Proof.}
Assume $\eta'<\eta$.

Since $\nu_{\omega^{\eta'}+1}$ is the largest element of $\overline{J}_{\omega^{\eta'}}$,
$\varphi_{\omega^{\eta'}}(\nu_{\omega^{\eta'}+1})$ is the largest element of the domain of $\iota_{\omega^{\eta'}}$. 
Since $\nu_{\omega^{\eta'}+1}$ is divisible by $\kappa_\alpha$, there exists $\delta$ such that $\nu_{\omega^{\eta'}+1}=\nu_{\omega^{\eta'}}+\kappa_\alpha\cdot\delta$.
Therefore, $\varphi_{\omega^{\eta'}}(\nu_{\omega^{\eta'}}+\kappa_\alpha\cdot\delta) = \alpha\cdot\delta$ is the largest element of the domain of $\iota_{\omega^{\eta'}}$.
Since $J_{\omega^{\eta'}}$ is $\rho$-isomorphic to a proper initial segment of $J_{\omega_\eta}$ by Lemma \ref{lorB}, $J_{\omega^{\eta'}}$ is $\rho$-isomorphic to $[\nu_{\omega^\eta},\nu_{\omega^\eta}+\kappa_\alpha\cdot\delta)$ and $\nu_{\omega^\eta}+\kappa_\alpha\cdot\delta\in J_{\omega^\eta}$.
Therefore, $\alpha\cdot\delta=\varphi_{\omega^\eta}(\nu_{\omega^\eta}+\kappa_\alpha\cdot\delta)$ is in the domain of $\iota_{\omega^\eta}$.
Since $\alpha\cdot\delta$ is the largest element of the domain of $\iota_{\omega^{\eta'}}$, $dom(\iota_{\omega^{\eta'}})\subseteq dom(\iota_{\omega^\eta})$.
\qed


\begin{lem} 
\label{lorC}
Assume $\xi<\theta_2$ and $\nu_\xi<^\rho_2 \nu$.
If $Y\subseteq [\nu_\xi,\nu)$ is $\rho$-closed and $\nu_\xi\in Y$ then there exists $\eta$ such that $\omega^\eta<\xi$  and a $\rho$-covering $h$ of $Y$ into $J_{\omega^\eta}$ such that $h(\nu_\xi)=\nu_{\omega^\eta}$ and $h(\mu)\leq^\rho_1 \nu_{\omega^\eta+1}$ whenever $\mu\in Y$ and $\mu\leq^\rho_1 \nu$.
\end{lem}
{\bf Proof.}
Assume $Y\subseteq [\nu_\xi,\nu)$ is $\rho$-closed and $\nu_\xi\in Y$. 

Since $\nu_\xi<^\rho_2 \nu$, parts 1 and 11 of Lemma 8.6 of [\ref{Ca??}] imply $\xi$ is a limit ordinal.
Therefore, there is a $\rho$-covering $h_1$ of $Y$ into $\nu_\xi$ such that $\nu_1<h_1[Y]$ and $h_1(\mu)\leq^\rho_1 \nu_\xi$ whenever $\mu\in Y$ and $\mu \leq^\rho_1 \nu$.

There exists $\xi'<\xi$ such that $h_1(\nu_\xi)\in J_{\xi'}$.
Since $\nu_1<h_1[Y]$, $1\leq \xi'$.

Define a function  $h_2$ on the range of $h_1$  such that 
\begin{equation}
h_2(\beta +\epsilon) =
\begin{cases}
\nu_{\xi'} +\epsilon & \text{if $\beta=h_1(\nu_\xi)$} \\
\beta+\epsilon & \text{if $\beta\not= h_1(\nu_\xi)$}
\end{cases}
\end{equation}
whenever $\beta+\epsilon$ is in the range of $h_1$, $\beta$ is divisible by $\rho$ and $\epsilon<\rho$.
By  Lemma 2.4 of [\ref{Ca??}], $h_2$ is an embedding of the range of $h_1$ into ${\cal R}^\rho$.
Clearly, the range of $h_2$ is contained in $[\nu_{\xi'},\nu_\xi)$.

Since $\nu_\xi\leq^\rho_1 \nu$, $h_1(\nu_\xi)\leq^\rho_1 \nu_\xi$.
By RTSI, $h_1(\nu_\xi)$ is divisible by $\kappa_\alpha$.
Therefore, $\nu_{\xi'}\leq^\rho_2 h_1(\nu_\xi)$.
This easily implies that $h_2$ is a $\rho$-covering of the range of $h_1$. 

Let $X$ be the intersection of $J_{\xi'}$ with the range of $h_2$.
Since $\nu_{\xi'+1}\leq^\rho_1 \nu_\xi$, there is a $\rho$-covering $h_3$ of the range of $h_2$ into $\nu_{\xi'+1}$ which fixes $X$. 

Since $\nu_{\xi'}$ is the least element of $X$ and $h_3(\nu_{\xi'})=\nu_{\xi'}$, $\nu_{\xi'}$ is the least element of the range of $h_3$.
Therefore, the range of $h_3$ is contained in $[\nu_{\xi'},\nu_{\xi'+1})=J_{\xi'}$.

Since $\xi'\not=0$, $\xi'=\zeta+\omega^\eta$ for some $\zeta$ and $\eta$.
By the First Recurrence Theorem for \leqtwo, there is an isomorphism $h_4$ of $\overline{J}_{\xi'}$ and $\overline{J}_{\omega^\eta}$. 

Let $h=h_4\circ h_3\circ h_2 \circ h_1$.

Since the composition of $\rho$-coverings is a $\rho$-covering, $h$ is a $\rho$-covering of $Y$.

A simple computation shows that $h(\nu_\xi)=\nu_{\omega^\eta}$.

Since the range of $h_3$ is contained in $J_{\xi'}$,  the range of $h$ is contained in $J_{\omega^\eta}$.

Assume $\mu\in Y$ and $\mu\leq^\rho_1 \nu$.
We will show $h(\mu)\leq^\rho_1 \nu_{\omega^\eta+1}$.
In the case $\mu=\nu_\xi$, $h(\mu)=\nu_{\omega^\eta}$ implying $h(\mu)\leq^\rho_1 \nu_{\omega^\eta+1}$.
So, we may assume $\nu_\xi<\mu$.

By choice of $h_1$, $h_1(\mu)\leq^\rho_1 \nu_\xi$.
By RTSI, $h_1(\mu)$ is divisible by $\kappa_\alpha$.
Since $\nu_\xi\leq^\rho_2 \nu$ and $\mu\leq^\rho_1\nu$, $\nu_\xi \leq^\rho_2 \mu$.
Therefore, $h_1(\nu_\xi)\leq^\rho_2 h_1(\mu)$.
By part 4 of Lemma 8.6 of [\ref{Ca??}], $h_1(\mu)<\nu_{\xi'+1}$ implying $h_1(\mu)\in J_{\xi'}$.
Since $h_1(\nu_\xi)<h_1(\mu)$, this implies $h_2(h_1(\mu))=h_1(\mu)$. 
Since $h_3$ fixes the elements of $X$, $h_3(h_1(\mu))=h_1(\mu)$.
Therefore, $h(\mu)=h_4(h_1(\mu))$.
Since $h_1(\mu)\leq^\rho_1\nu_\xi$, $h_1(\mu)\leq^\rho_1 \nu_{\xi'+1}$.
Since $h_4$ is an isomorphism of $\overline{J}_{\xi'}$ and $\overline{J}_{\omega^\eta}$, $h(\mu)=h_4(h_1(\mu))\leq^\rho_1 h_4(\nu_{\xi'+1})=\nu_{\omega^\eta+1}$.
\qed

\blankline

The Second Recurrence Theorem for \leqtwo\ below provides a computation of the lengths of the intervals $\overline{J}_\xi$ in terms of the  intervals $I^\alpha_\beta$.
By the First Recurrence Theorem for \leqtwo, this reduces to calculating the lengths of the intervals $\overline{J}_\xi$ when $\xi$ has the form $\omega^\eta$.


\begin{thm}
{\bf (Second Recurrence Theorem for \leqtwo)} 
 If $\omega^\eta<\theta_2$ and $\eta=\lambda+n$ where $n\in \omega$ and $\lambda$ is a multiple of $\omega$, possibly 0, then $\alpha\cdot\lambda<\theta^\alpha_1$ and the following hold.
\begin{enumerate}
\item
If $\omega^\eta+1<\theta_2$ then the following hold.
\begin{enumerate}
\item
If the largest element of $I^\alpha_{\alpha\cdot\lambda}$ is not a successor multiple of $\alpha$ then the largest element of the domain of $\iota_{\omega^\eta}$ is $max(I^\alpha_{\alpha\cdot(\eta+1)})$ i.e. $\kappa^\alpha_{\alpha\cdot(\eta+1)}$.
\item
If the largest element of $I^\alpha_{\alpha\cdot\lambda}$ is  a successor multiple of $\alpha$ then the largest element of the domain of $\iota_{\omega^\eta}$ is $max(I^\alpha_{\alpha\cdot\eta})$ which is $\kappa^\alpha_{\alpha\cdot\eta}$ if $n>0$.
 \end{enumerate}
 \item
If $\omega^\eta+1=\theta_2$ and $\alpha+1<\theta_1$ then the following hold.
 \begin{enumerate}
\item
If the largest element of $I^\alpha_{\alpha\cdot\lambda}$ is not a successor multiple of $\alpha$ then the largest element of the domain of $\iota_{\omega^\eta}$ is $max(I^\alpha_{\alpha\cdot\eta+\gamma})$ for some $\gamma<\alpha$.
\item
If the largest element of $I^\alpha_{\alpha\cdot\lambda}$ is  a successor multiple of $\alpha$ then $\eta$ is a successor ordinal and the largest element of the domain of $\iota_{\omega^\eta}$ is $max(I^\alpha_{\alpha\cdot\eta'+\gamma})$ for some $\gamma<\alpha$ where $\eta=\eta'+1$.
 \end{enumerate}
 \item If $\omega^\eta+1=\theta_2$ and $\alpha+1=\theta_1$ then  $\eta+1=\theta^\alpha_1$.
\end{enumerate}
\end{thm}
{\bf Proof.}
Recall that Theorem 8.8 of [\ref{Ca??}] says that either $\theta_2=\infty$ or $\theta_2=\theta+1$ where $\theta$ is infinite and additively indecomposable.

Notice that for $\xi<\theta_2$, if $J_\xi$ is unbounded iff $\alpha+1=\theta_1$ and  $\xi+1=\theta_2$. Moreover, $J_\xi$ is unbounded iff the domain of $\iota_\xi$ is unbounded.

We claim that if $\xi+1<\theta_2$ then the largest element of the domain of $\iota_\xi$ is a successor multiple of $\alpha$.
To see this, suppose $\xi+1<\theta_2$. 
The largest element of $\overline{J}_\xi$ is $\nu_{\xi+1}$ which is a successor multiple of $\kappa_\alpha$ by part 11 of Lemma 8.6 of [\ref{Ca??}].
Therefore, the largest element of the domain of $\iota_\xi$, $\varphi_\xi(\nu_{\xi+1})$, is a successor multiple of $\alpha$.

When $J_{\omega^\eta}$ is bounded, define $f(\eta)$ to be the largest $\beta$ such that $\kappa^\alpha_\beta\in dom(\iota_{\omega^\eta})$.
By part 3 of Lemma \ref{lorA}, if $\xi<\theta_2$ and the domain of $\iota_\xi$ intersects $I^\alpha_\beta$ then $dom(\iota_\xi)$ contains $I^\alpha_\beta$.
Therefore, if $J_{\omega^\eta}$ is bounded then the largest element of $dom(\iota_{\omega^\eta})$ is $max(I^\alpha_{f(\eta)})$.
By the previous paragraph, if $\omega^\eta+1<\theta_2$ then the largest element of $I^\alpha_{f(\eta)}$ is a successor multiple of $\alpha$.

We will now establish the theorem by induction on $\eta$. Let $K$ be the class of ordinals $\eta$ such that the theorem holds.

Assume $\eta$ is an ordinal such that $\eta'\in K$ whenever $\eta'<\eta$ and suppose $\omega^\eta<\theta_2$ and $\eta=\lambda+n$ where $n\in \omega$ and $\lambda$ is a multiple of $\omega$.

\halfblankline

{\bf Case 1.} 
Assume $\eta=0$.

In this case, $\lambda=n=0$ and $\omega^\eta=1$.

Since $I^\alpha_{\alpha\cdot \lambda}=I^\alpha_0=\{0\}$ by part 1 of Lemma 5.4 of [\ref{Ca??}], and $0$ is not a successor multiple of $\alpha$, the only clause which is not vacuously true is 2(a). 

By parts 1 and 11 of Lemma 8.6 of [\ref{Ca??}], $\nu_m=\kappa_\alpha\cdot(m+1)$ for $m\in \omega$.
Therefore, $\overline{J}_{\omega^\eta}=\overline{J}_1=[\nu_1,\nu_2]=[\kappa_\alpha\cdot 2, \kappa_\alpha\cdot 3]$ implying that the domain of $\iota_{\omega^\eta}$ is $[0,\alpha]$.
By parts 1 and 8 of Lemma 8.6 of [\ref{Ca??}], $I^\alpha_{\alpha\cdot(\eta+1)}=I^\alpha_\alpha =\{\alpha\}$.
Therefore, the largest element of the domain of $\iota_{\omega^\eta}$ is $max(I_{\alpha\cdot(\eta+1)})$.

\halfblankline

{\bf Case 2.}
Assume $\eta$ is a successor ordinal.

Let $\eta=\eta'+1$. 
In this case, $n>0$ and $\eta'=\lambda+(n-1)$.
By the induction hypothesis, $\alpha\cdot\lambda<\theta^\alpha_1$.

Since $\omega^{\eta'}+1<\theta_2$, the largest element of $\overline{J}_{\omega^{\eta'}}$ is $\nu_{\omega^{\eta'}+1}$ which is a successor multiple of $\kappa_\alpha$.
Therefore, there exists a successor ordinal $\delta$ such that the largest element of $\overline{J}_{\omega^{\eta'}}$ is $\nu_{\omega^{\eta'}}+\kappa_\alpha\cdot\delta$.
Therefore, $max(dom(\iota_{\omega^{\eta'}}))=\varphi_{\omega^{\eta'}}(\nu_{\omega^{\eta'}+1})=\alpha\cdot\delta$. 
Since $\overline{J}_{\omega^{\eta'}}$ is bounded, $f(\eta')$ is defined and $max(dom(\iota_{\omega^{\eta'}}))=max(I^\alpha_{f(\eta')})$.
Therefore, $max(I^\alpha_{f(\eta')})=\alpha\cdot\delta$.

\halfblankline

{\bf Claim 1 for Case 2.}
$J_{\omega^{\eta'}}\cong_\rho [\nu_{\omega^\eta},\nu_{\omega^\eta}+\kappa_\alpha\cdot\delta)$ and $\nu_{\omega^\eta}+\kappa_\alpha\cdot\delta\in J_{\omega^\eta}$.

By Lemma \ref{lorB}, $J_{\omega^{\eta'}}$ is $\rho$-isomorphic to a proper initial segment of $J_{\omega^\eta}$.
Since $J_{\omega^{\eta'}}=[\nu_{\omega^{\eta'}}, \nu_{\omega^{\eta'}}+\kappa_\alpha\cdot\delta)$, this implies that $J_{\omega^{\eta'}}\cong_\rho [\nu_{\omega^\eta},\nu_{\omega^\eta}+\kappa_\alpha\cdot\delta)$ and $\nu_{\omega^\eta}+\kappa_\alpha\cdot\delta\in J_{\omega^\eta}$. 

\halfblankline

By Lemma 6.7 of [\ref{Ca??}], there is an $\alpha$-incompressible set $Z$ whose largest element is $\alpha\cdot\delta$.

Since $\varphi^{-1}_{\omega^\eta}[[0,\alpha\cdot\delta]]=[\nu_{\omega^\eta},\nu_{\omega^\eta}+\kappa_\alpha\cdot\delta]$, part 2 of Lemma \ref{lsrtt} implies there exists $Z^*\subseteq \varphi^{-1}_{\omega^\eta}[[0,\alpha\cdot\delta]]$ which is locally incompressible such that $\nu_{\omega^\eta}\in Z^*$ and $\iota_{\omega^\eta}[Z]\subseteq Z^*$. 

\halfblankline

{\bf Claim 2 for Case 2.}
$\nu_{\omega^\eta}+\kappa_\alpha\cdot \delta + \kappa_\alpha\not\in J_{\omega^\eta}$.

To prove the claim, argue by contradiction and assume $\nu_{\omega^\eta}+\kappa_\alpha\cdot \delta + \kappa_\alpha\in J_{\omega^\eta}$.

Since $\nu_{\omega^\eta}+\kappa_\alpha\cdot\delta+\kappa_\alpha$ is divisible by $\kappa_\alpha$, $\nu_{\omega^\eta}\leq^\rho_2 \nu_{\omega^\eta}+\kappa_\alpha\cdot\delta+\kappa_\alpha$. 
By Lemma \ref{lorC}, there is a $\zeta<\eta$ and $\rho$-covering of $Z^*$ into $J_{\omega^\zeta}$  which maps $\nu_{\omega^\eta}$ to $\nu_{\omega^\zeta}$. 
By Lemma \ref{lorB}, $J_{\omega^\zeta}$ is $\rho$-isomorphic to an initial segment of $J_{\omega^{\eta'}}$.
Therefore, we may assume $\zeta=\eta'$.
Since $J_{\omega^{\eta'}}\cong_\rho [\nu_{\omega^\eta},\nu_{\omega^\eta}+\kappa_\alpha\cdot\delta)$, there is a $\rho$-covering $h$ of $Z^*$ into $[\nu_{\omega^\eta},\nu_{\omega^\eta}+\kappa_\alpha\cdot\delta)$ such that $h(\nu_{\omega^\eta})=\nu_{\omega^\eta}$.

By part 2 of Lemma \ref{lorA}, there is an $\alpha$-covering $h'$ of $Z$ into $[0,\alpha\cdot\delta]$ such that $h'(\sigma+\gamma)= (\iota_{\omega^\eta}^{-1}\circ h \circ \iota_{\omega^\eta})(\sigma)+\gamma$ whenever $\sigma+\gamma\in Z$, $\sigma$ is divisible by $\alpha$ and $\gamma<\alpha$.

Since $Z$ in $\alpha$-incompressible and $\alpha\cdot\delta$ is the largest element of $Z$, we we will contradict part 9 of Lemma 6.4 of [\ref{Ca??}] by showing $h'(\alpha\cdot\delta)<\alpha\cdot\delta$. 
This inequality follows from the definition of $h'$ and the fact that  $h(\iota_\xi(\alpha\cdot\delta))\in h[Z^*]<\nu_{\omega^\eta}+\kappa_\alpha\cdot\delta$.

\halfblankline

To verify part 1, assume $\omega^\eta+1<\theta_2$.

The assumption $\omega^\eta+1<\theta_2$ implies $J_{\omega^\eta}=[\nu_{\omega^\eta},\nu_{\omega^\eta+1})$. 
Since $\nu_{\omega^\eta+1}$ is a successor multiple of $\kappa_\alpha$, Claims 1 and 2 imply that $\nu_{\omega^\eta+1}=\nu_{\omega^\eta}+\kappa_\alpha\cdot\delta + \kappa_\alpha$.
Therefore, 
\begin{tabbing}
\hspace{1.2in} $max(dom(\iota_{\omega^\eta}))$ \ \= = \= $\varphi_{\omega^\eta}(\nu_{\omega^\eta+1})$ \\
\> = \> $\varphi_{\omega^\eta}(\nu_{\omega^\eta}+\kappa_\alpha\cdot\delta+\kappa_\alpha)$ \\
\> = \> $\alpha\cdot\delta+\alpha$ \\
\> = \> $max(I^\alpha_{f(\eta')})+\alpha$  \\
\> = \> $max(I^\alpha_{f(\eta')+\alpha})$ 
\end{tabbing}

To establish part 1(a), assume $max(I^\alpha_{\alpha\cdot\lambda})$ is not a successor multiple of $\alpha$. 
By the induction hypothesis, $f(\eta')=\alpha\cdot(\eta'+1)=\alpha\cdot\eta$.
By the string of equalities above, the largest element of the domain of $\iota_{\omega^\eta}$ is $max(I^\alpha_{\alpha\cdot(\eta+1)})$ as desired.

To establish part 1(b), assume $max(I^\alpha_{\alpha\cdot\lambda})$ is a successor multiple of $\alpha$.
By the induction hypothesis, $f(\eta')=\alpha\cdot\eta'$.
By the string of equalities above, the largest element of the domain of $\iota_{\omega^\eta}$ is $max(I^\alpha_{\alpha\cdot\eta})$ as desired.

This concludes the proof of part 1 for Case 2. 

\halfblankline

To verify part 2, assume $\omega^\eta+1=\theta_2$.

The assumption $\omega^\eta+1=\theta_2$ implies $\overline{J}_{\omega^\eta}=J_{\omega^\eta}$.
Using Claim 2, this implies that $\alpha\cdot\delta + \alpha \not\in dom(\iota_{\omega^\eta})$.
On the other hand, Claim 1 implies that $\alpha\cdot\delta\in dom(\iota_{\omega^\eta})$.
Therefore, there exists $\gamma<\alpha$ such that $max(dom(\iota_{\omega^\eta}))=\alpha\cdot\delta+\gamma$.
Similar to the computation above, we have
\begin{tabbing}
\hspace{1.2in} $max(dom(\iota_{\omega^\eta}))$ \ \= = \= $\alpha\cdot\delta+\gamma$ \\
\> = \> $max(I^\alpha_{f(\eta')})+\gamma$  \\
\> = \> $max(I^\alpha_{f(\eta')+\gamma})$ 
\end{tabbing}

To establish part 2(a), assume $max(I^\alpha_{\alpha\cdot\lambda})$ is not a successor multiple of $\alpha$. 
By the induction hypothesis, $f(\eta')=\alpha\cdot(\eta'+1)=\alpha\cdot\eta$.
By the string of equalities above, the largest element of the domain of $\iota_{\omega^\eta}$ is $max(I^\alpha_{\alpha\cdot\alpha+\gamma})$ and $\gamma<\alpha$ as desired.

To establish part 2(b), assume $max(I^\alpha_{\alpha\cdot\lambda})$ is a successor multiple of $\alpha$.
By the induction hypothesis, $f(\eta')=\alpha\cdot\eta'$.
By the string of equalities above, the largest element of the domain of $\iota_{\omega^\eta}$ is $max(I^\alpha_{\alpha\cdot\eta'+\gamma})$ and $\gamma<\alpha$ as desired.

This concludes the proof of part 2 for Case 2.

Since Claim 2 implies that $J_{\omega^\eta}$ is bounded, the hypothesis of part 3 fails and part 3 holds vacuously.

\halfblankline

{\bf Case 3.}
Assume $\eta$ is a limit ordinal.

In this case, $n=0$ and $\lambda=\eta$.

By the induction hypothesis, part 1 of the theorem holds with $\eta$ replaced by $\eta'$ for all $\eta'<\eta$.
Since $f(\eta')$ is either $\alpha\cdot\eta'$ or $\alpha\cdot(\eta'+1)$ for all $\eta'<\eta$, $\alpha\cdot\eta$ is the least upper bound of the ordinals $f(\eta')$ for $\eta'<\eta$ and $\alpha\cdot\eta<\theta^\alpha_1$.
Moreover, if $\eta_1<\eta_2<\eta$ then $f(\eta_1)<f(\eta_2)$.

By Lemma \ref{lorB1}, $\kappa^\alpha_{f(\eta')}\in dom(\iota_{\omega^\eta})$ for $\eta'<\eta$. 
Since the domain of $\iota_\eta$ is a closed interval, this implies that $\kappa^\alpha_{\alpha\cdot\eta}\in dom(\iota_{\omega^\eta})$.
Therefore, $I^\alpha_{\alpha\cdot\eta}$ is a subset of $dom(\iota_{\omega^\eta})$.

First, suppose $I^\alpha_{\alpha\cdot\eta}$ is unbounded.
This implies that $\alpha\cdot\eta+1=\theta^\alpha_1$.
By Corollary 7.3 of [\ref{Ca??}], $\alpha\cdot\eta$ is an epsilon number greater than $\alpha$ implying $\alpha\cdot\eta=\eta$ and, hence, $\eta+1=\theta^\alpha_1$.
Since $I^\alpha_{\alpha\cdot\eta}$ a subset of the domain of $\iota_{\omega^\eta}$, $J_{\omega^\eta}$ is unbounded implying $\omega^\eta+1=\theta_2$.
Therefore, the conclusion of part 3 holds. 
Parts 1 and 2 hold vacuously.

For the remainder of the proof, assume $I^\alpha_{\alpha\cdot\eta}$ is bounded.
This makes part 3 vacuously true.

Write the largest element of $I^\alpha_{\alpha\cdot\eta}$ as  $\alpha\cdot\delta+\gamma$ where $\gamma<\alpha$.

\halfblankline

{\bf Claim 1 for Case 3.}
$\nu_{\omega^\eta}+\kappa_\alpha\cdot\delta+\kappa_\gamma\in \overline{J}_{\omega^\eta}$.

This follows from the fact $\nu_{\omega^\eta}+\kappa_\alpha\cdot\delta + \kappa_\gamma = \iota_{\omega^\eta}(\alpha\cdot\delta+\gamma)$.

\halfblankline

{\bf Claim 2 for Case 3.}
If $\alpha\cdot\delta+\gamma$ is not a successor multiple of $\alpha$ then $\nu_{\omega^\eta}+\kappa_\alpha\cdot\delta + \kappa_\gamma\in J_{\omega^\eta}$.

Assume $\alpha\cdot\delta+\gamma$ is not a successor multiple of $\alpha$.
This implies $\nu_{\omega^\eta}+\kappa_\alpha\cdot\delta+\kappa_\gamma$ is not a successor multiple of $\kappa_\alpha$.
By Claim 1, $\nu_{\omega^\eta}+\kappa_\alpha\cdot\delta +\kappa_\gamma\in \overline{J}_{\omega^\eta}$.
In case $\omega^\eta+1=\theta_2$, $\overline{J}_{\omega^\eta}=J_{\omega^\eta}$ and the desired conclusion follows. 
Suppose $\omega^\eta+1<\theta_2$. Since $\overline{J}_{\omega^\eta}=[\nu_{\omega^\eta}\nu_{\omega^\eta+1}]$ and $\nu_{\omega^\eta+1}$ is a successor multiple of $\kappa_\alpha$, $\nu_{\omega^\eta}+\kappa_\alpha\cdot\delta +\kappa_\gamma\in [\nu_{\omega^\eta},\nu_{\omega^\eta+1})=J_{\omega^\eta}$.

\halfblankline

By Lemma 6.7 of [\ref{Ca??}], there is an $\alpha$-incompressible set $Z$ such that $max(Z)=\alpha\cdot\delta+\gamma$. 
By part 2 of Lemma \ref{lsrtt}, there is a locally incompressible subset $Z^*$ of
$[\nu_{\omega^\eta},\nu_{\omega^\eta}+\kappa_\alpha\cdot\delta+\kappa_{\gamma+1})$  such that $\iota_{\omega^\eta}[Z]\subseteq Z^*$ and $\nu_{\omega^\eta}\in Z^*$.
By part 13 of Lemma 8.6 of [\ref{Ca??}], $[\nu_{\omega^\eta},\nu_{\omega^\eta}+\kappa_\alpha\cdot\delta+\kappa_{\gamma+1})\subseteq\overline{J}_{\omega^\eta}$.
Therefore, $Z^*$ is a subset of $\overline{J}_{\omega^\eta}$. 

\halfblankline

{\bf Claim 3 for Case 3.}
$\nu_{\omega^\eta}+\kappa_\alpha\cdot\delta+\kappa_\alpha\not\in J_{\omega^\eta}$.

The proof is similar to that of Claim 2 for Case 2.

Argue by contradiction and assume $\nu_{\omega^\eta}+\kappa_\alpha\cdot \delta + \kappa_\alpha\in J_{\omega^\eta}$.

Since $\nu_{\omega^\eta}+\kappa_\alpha\cdot\delta+\kappa_\alpha$ is divisible by $\kappa_\alpha$, $\nu_{\omega^\eta}\leq^\rho_2 \nu_{\omega^\eta}+\kappa_\alpha\cdot\delta+\kappa_\alpha$. 
By Lemma \ref{lorC}, there is a $\eta'<\eta$ and $\rho$-covering $h_1$ of $Z^*$ into $J_{\omega^{\eta'}}$  which maps $\nu_{\omega^\eta}$ to $\nu_{\omega^{\eta'}}$. 
Since $\nu_{\omega^{\eta'}+1}$ is a multiple of $\kappa_\alpha$, there exists $\delta'$ such that $\nu_{\omega^{\eta'}+1}=\nu_{\omega^\eta}+\kappa_\alpha\cdot\delta'$.
Since $max(\overline{J}_{\omega^{\eta'}})=\nu_{\omega^{\eta'}+1}$, this implies $max(I^\alpha_{f(\eta')})=max(dom(\iota_{\omega^{\eta'}}))=\varphi_{\omega^{\eta'}}(\nu_{\omega^\eta}+\kappa_\alpha\cdot\delta')=\alpha\cdot\delta'$.
Since $f(\eta')<\alpha\cdot\eta$, this implies $\alpha\cdot\delta'<\alpha\cdot\delta+\gamma$.

By Lemma \ref{lorB}, $J_{\omega^{\eta'}}$ is $\rho$-isomorphic to a proper initial segment of $J_{\omega^\eta}$.
Let $h_2$ be the $\rho$-isomorphism of  $J_{\omega^{\eta'}}$ and $ [\nu_{\omega^\eta},\nu_{\omega_\eta}+\kappa_\alpha\cdot\delta')$.
Let $h=h_2\circ h_1$.
$h$ is a $\rho$-covering of $Z^*$ into $[\nu_{\omega^\eta},\nu_{\omega^\eta}+\kappa_\alpha\cdot\delta')$ and $h(\nu_{\omega^\eta})=\nu_{\omega_\eta}$.
By part 2 of Lemma \ref{lorA}, there is an $\alpha$-covering $h'$ of $Z$ into $[0,\alpha\cdot\delta+\gamma]$ such that $h'(\sigma+\zeta)=(\iota_{\omega^\eta}^{-1}\circ h \circ \iota_{\omega^\eta})(\sigma)+\zeta$ whenever $\sigma +\zeta\in Z$, $\sigma$ is divisible by $\alpha$ and $\zeta<\alpha$.

We claim $h'$ maps $Z$ into $[0,\alpha\cdot\delta')$ which contradicts fact $Z$ is $\alpha$-incompressible. 
To see this, it suffices to show $h'(\alpha\cdot\delta+\gamma)<\alpha\cdot\delta'$. 
Since $h(\iota_\xi(\alpha\cdot\delta))\in h[Z*]<\nu_{\omega^{\omega^\eta}}+\kappa_\alpha\cdot\delta'$, the definition of $h'$ implies $h'(\alpha\cdot\delta)<\alpha\cdot\delta'$.
Since $\alpha$ is additively indecomposable, $h'(\alpha\cdot\delta+\gamma)=h'(\alpha\cdot\delta)+\gamma<\alpha\cdot\delta'$.

\halfblankline

{\bf Claim 4 for Case 3.}
If $\alpha\cdot\delta+\gamma$ is a successor multiple of $\alpha$ then $\nu_{\omega^\eta}+\kappa_\alpha\cdot\delta\not\in J_{\omega^\eta}$.

Argue by contradiction and assume $\alpha\cdot\delta+\gamma$ is a successor multiple of $\alpha$ but $\nu_{\omega^\eta}+\kappa_\alpha\cdot\delta\in J_{\omega^\eta}$.
This implies $\delta$ is a successor ordinal and $\gamma=0$.
Therefore, $Z^*\subseteq [\nu_{\omega^\eta},\nu_{\omega^\eta}+\kappa_\alpha\cdot\delta+\kappa_1)=[\nu_{\omega^\eta},\nu_{\omega^\eta}+\kappa_\alpha\cdot\delta]$. 

Since $\nu_{\omega^\eta}+\kappa_\alpha\cdot\delta$ is divisible by $\kappa_\alpha$,  $\nu_{\omega^\eta}\leq^\rho_2 \nu_{\omega^\eta}+\kappa_\alpha\cdot\delta$.
Let $Z^-$ be $Z^*-\{\nu_{\omega^\eta}+\kappa_\alpha\cdot\delta\}$.
Since $Z^*\subseteq [\nu_{\omega^\eta},\nu_{\omega^\eta}+\kappa_\alpha\cdot\delta]$, $Z^-\subseteq [\nu_{\omega^\eta},\nu_{\omega^\eta}+\kappa_\alpha\cdot\delta)$.
By Lemma \ref{lorC}, there is a $\eta'<\eta$ and $\rho$-covering $h_1$ of $Z^-$ into $J_{\omega^{\eta'}}$  such that  $h_1(\nu_{\omega^\eta})=\nu_{\omega^{\eta'}}$ and $h_1(\mu)\leq^\rho_1 \nu_{\omega^{\eta'}+1}$ whenever $\mu\in Z^-$ and $\mu\leq^\rho_1 \nu_{\omega^\eta}+\kappa_\alpha\cdot\delta$ . 
Since $\nu_{\omega^{\eta'}+1}$ is a multiple of $\kappa_\alpha$, there exists $\delta'$ such that $\nu_{\omega^{\eta'}+1}=\nu_{\omega^\eta}+\kappa_\alpha\cdot\delta'$.
Since $max(\overline{J}_{\omega^{\eta'}})=\nu_{\omega^{\eta'}+1}$, this implies $max(I^\alpha_{f(\eta')})=max(dom(\iota_{\omega^{\eta'}}))=\varphi_{\omega^{\eta'}}(\nu_{\omega^\eta}+\kappa_\alpha\cdot\delta')=\alpha\cdot\delta'$.
Since $f(\eta')<\alpha\cdot\eta$, this implies $\alpha\cdot\delta'<\alpha\cdot\delta$.

By Lemma \ref{lorB}, $J_{\omega^{\eta'}}$ is $\rho$-isomorphic to a proper initial segment of $J_{\omega^\eta}$.
Let $h_2$ be the $\rho$-isomorphism of  $J_{\omega^{\eta'}}$ and $ [\nu_{\omega^\eta},\nu_{\omega_\eta}+\kappa_\alpha\cdot\delta')$.
Define $h$ on $Z^*$ by
\begin{equation}
h(\beta) =
\begin{cases}
h_2(h_1(\beta)) & \text{if $\beta\not=\nu_{\omega^\eta}+\kappa_\alpha\cdot\delta$} \\
\nu_{\omega^\eta}+\kappa_\alpha\cdot\delta'  & \text{if $\beta=\nu_{\omega^\eta}+\kappa_\alpha\cdot\delta$}
\end{cases}
\end{equation}
Notice that the range of $h$ is contained in $[\nu_{\omega^\eta},\nu_{\omega^\eta}+\kappa_\alpha\cdot\delta']$.

We claim that $h$ is a $\rho$-covering of $Z^*$.
By  Lemma 2.4 of [\ref{Ca??}], $h$ is an embedding  of $Z^*$ into ${\cal R}^\rho$.

Clearly, the restriction of $h$ to $Z^-$ is a $\rho$-covering of $Z^-$.

Since $\delta$ is a successor ordinal, $\nu_{\omega^\eta}+\kappa_\alpha\cdot\delta$ is not a limit of ordinals which are divisible by $\kappa_\alpha$.
Since $\nu_{\omega^\eta}\leq^\rho_2 \mu$ implies $\mu$ is divisible by $\kappa_\alpha$, $\nu_{\omega^\eta}+\kappa_\alpha\cdot\delta$ is not a limit of ordinals $\mu$ with $\nu_{\omega^\eta}\leq^\rho_2 \mu$.
By part 2 of Lemma 3.6 of [\ref{Ca??}], there is no $\mu$ such that $\nu_{\omega^\eta}<\mu$ and $\mu<^\rho_2 \nu_{\omega^\eta}+\kappa_\alpha\cdot\delta$.

Since $\nu_{\omega^\eta}+\kappa_\alpha\cdot\delta'<\nu_{\omega^\eta}+\kappa_\alpha\cdot\delta\in \overline{J}_{\omega^\eta}$, $\nu_{\omega^\eta}+\kappa_\alpha\cdot\delta'\in J_{\omega^\eta}$ and $\nu_{\omega^\eta}\leq^\rho_2 \nu_{\omega^\eta}+\kappa_\alpha\cdot\delta'$.
This fact and the conclusions of the previous two paragraphs easily imply $h$ is a $\rho$-covering of $Z^*$.

By part 2 of Lemma \ref{lorA}, there is an $\alpha$-covering $h'$ of $Z$ into $[0,\alpha\cdot\delta]$ such that $h'(\sigma+\zeta)=(\iota_{\omega^\eta}^{-1}\circ h \circ \iota_{\omega^\eta})(\sigma)+\zeta$ whenever $\sigma +\zeta\in Z$, $\sigma$ is divisible by $\alpha$ and $\zeta<\alpha$.

Clearly, $h'(\alpha\cdot\delta)=\alpha\cdot\delta'$.
Since $\alpha\cdot\delta=max(Z)$ and $\delta'<\delta$, this contradicts the fact that $Z$ is $\alpha$-incompressible.

\halfblankline

To verify part 1, assume $\omega^\eta+1<\theta_2$.
This implies that $J_{\omega^\eta}=[\nu_{\omega^\eta},\nu_{\omega^\eta+1})$. 

To establish part 1(a), assume $\alpha\cdot\delta+\gamma$ is not a successor multiple of $\alpha$. 
We will show that $max(dom(\iota_{\omega^\eta}))=max(I^\alpha_{\alpha\cdot(\eta+1)})$.
Since $\nu_{\omega^\eta+1}$ is a successor multiple of $\kappa_\alpha$, Claims 2 and 3 imply that $\nu_{\omega^\eta+1}=\nu_{\omega^\eta}+\kappa_\alpha\cdot\delta + \kappa_\alpha$.
Therefore, 
\begin{tabbing}
\hspace{1.2in} $max(dom(\iota_{\omega^\eta}))$ \ \= = \= $\varphi_{\omega^\eta}(\nu_{\omega^\eta+1})$ \\
\> = \> $\varphi_{\omega^\eta}(\nu_{\omega^\eta}+\kappa_\alpha\cdot\delta+\kappa_\alpha)$ \\
\> = \> $\alpha\cdot\delta+\alpha$ \\
\>=\> $\alpha\cdot\delta+\gamma+\alpha$ \\
\> = \> $max(I^\alpha_{\alpha\cdot\eta})+\alpha$ \\
\> = \> $max(I^\alpha_{\alpha\cdot\eta +\alpha})$  
\end{tabbing}

To establish part 1(b), assume $\alpha\cdot\delta+\gamma$ is a successor multiple of $\alpha$.
We will show that $max(dom(\iota_{\omega^\eta}))=max(I^\alpha_{\alpha\cdot\eta})$.
In this case $\gamma=0$.
Since $\nu_{\omega^\eta+1}$ is a successor multiple of $\kappa_\alpha$, Claims 1 and 4 imply that $\nu_{\omega^\eta+1}=\nu_{\omega^\eta}+\kappa_\alpha\cdot\delta$.
Therefore, 
\begin{tabbing}
\hspace{1.2in} $max(dom(\iota_{\omega^\eta}))$ \ \= = \= $\varphi_{\omega^\eta}(\nu_{\omega^\eta+1})$ \\
\> = \> $\varphi_{\omega^\eta}(\nu_{\omega^\eta}+\kappa_\alpha\cdot\delta)$ \\
\> = \> $\alpha\cdot\delta$ \\
\>=\> $\alpha\cdot\delta$ \\
\> = \> $max(I^\alpha_{\alpha\cdot\eta})$   
\end{tabbing}

This concludes the proof of part 1 for Case 2. 

\halfblankline

To verify part 2, assume $\omega^\eta+1=\theta_2$.

The assumption $\omega^\eta+1=\theta_2$ implies $\overline{J}_{\omega^\eta}=J_{\omega^\eta}$.
By Claims 1 and 4, $\alpha\cdot\delta+\gamma$ is not a successor multiple of $\alpha$ making part 2(b) vacuously true.

By Claim 3 and the fact $\overline{J}_{\omega^\eta}=J_{\omega^\eta}$, $\alpha\cdot\delta+\alpha\not\in dom(\iota_{\omega^\eta})$.
Since $\alpha\cdot\delta+\gamma\in dom(\iota_{\omega^\eta})$, this implies there exists $\gamma'<\alpha$ such that $max(dom(\iota_{\omega^\eta}))=\alpha\cdot\delta+\gamma+\gamma'$. 
To establish part 2(a), we will show that $max(dom(\iota_{\omega^\eta}))=max(I^\alpha_{\alpha\cdot\delta +\gamma'})$.
\begin{tabbing}
\hspace{1.2in} $max(dom(\iota_{\omega^\eta}))$ \ \= = \= $\alpha\cdot\delta+\gamma+\gamma'$ \\
\> = \> $max(I^\alpha_{\alpha\cdot\eta})+\gamma'$  \\
\> = \> $max(I^\alpha_{\alpha\cdot\eta+\gamma'})$ 
\end{tabbing}
\qed

\vspace{10 mm}


\begin{center}
 REFERENCES
\end{center}

\begin{enumerate}
\item
\label{Ba75}
 J. Barwise, {\bf Admissible Sets and Structures}, 
Springer-Verlag, Berlin, 1975.
\item
\label{Ca09}
T. Carlson, {\it Patterns of resemblance of order 2}, Annals of Pure and Applied Logic 158 (2009), pp. 90-124.
\item
\label{Ca16}
T. Carlson, {\it Generalizing Kruskal's Theorem to pairs of cohabitating trees},Archive for Mathematical Logic 55 (2016), pp. 37-48.
\item
\label{Ca??}
T. Carlson, {\it Structural Properties of \calRtwo\ Part I}, preprint.
\end{enumerate}

\end{document}